\documentclass[11pt,reqno]{amsart}
\usepackage{amsfonts}
\usepackage{amsmath,amssymb,amsthm}\allowdisplaybreaks[4]
\usepackage{url}
\usepackage{graphicx}   
\usepackage{verbatim}       
\usepackage{layout}
\usepackage{galois}
\usepackage{mathrsfs}
\usepackage{bbm}
\usepackage{color}
\usepackage{cite}
\numberwithin{equation}{section}
\usepackage[linkcolor=red,citecolor=blue]{hyperref}
\usepackage{enumerate}

\newtheorem{theorem}{Theorem}[section]
\newtheorem{lemma}[theorem]{Lemma}
\newtheorem{proposition}[theorem]{Proposition}
\newtheorem{remark}[theorem]{Remark}
\newtheorem{definition}[theorem]{Definition}
\newtheorem{corollary}[theorem]{Corollary}

\newtheorem{example}[theorem]{Example}

\newtheorem{assumption}[theorem]{Assumption}

\def\[{\left[}\def\]{\right]}
\def\<{\left<}\def\>{\right>}
\def\({\left(}\def\){\right)}

\def\beqlb{\begin{eqnarray}}\def\eeqlb{\end{eqnarray}}
\def\beqnn{\begin{eqnarray*}}\def\eeqnn{\end{eqnarray*}}

\def\ar{\!\!&}

\def\al{{\alpha}}\def\be{{\beta}}
\def\ep{{\epsilon}}\def\ga{{\gamma}}\def\ka{{\kappa}}

\def\La{{\Lambda}}

\begin{document}
\title{Speed of coming down from infinity for $\Lambda$-Fleming-Viot initial support}
\author{ Huili Liu and Xiaowen Zhou}
\address{Huili Liu: School of Mathematical Sciences, Hebei Normal University,
Shijiazhuang, Hebei, 050024, China}
\email{liuhuili@hebtu.edu.cn}
\address{Xiaowen Zhou: Department of Mathematics and Statistics, Concordia University,
1455 de Maisonneuve Blvd. West, Montreal, Quebec, H3G 1M8, Canada}
\email{xiaowen.zhou@concordia.ca }
\date{\today}

\begin{abstract}
The $\Lambda$-Fleming-Viot process is a probability measure-valued process that is dual to a $\Lambda$-coalescent that allows multiple collisions. In this paper, we consider a class of $\Lambda$-Fleming-Viot processes with Brownian spatial motion and with associated
$\Lambda$-coalescents that come down from infinity. Notably, these processes have the compact support property: the support of the process becomes finite as soon as $t>0$, even though the initial measure has unbounded support.
We obtain asymptotic results characterizing the rates at which the initial supports become finite. The  rates of coming down are expressed in terms of the asymptotic inverse function of the tail distribution of the initial measure and the speed function of coming down from infinity for the corresponding $\Lambda$-coalescent.

\end{abstract}
\maketitle \pagestyle{myheadings} \markboth{\textsc{Speed of CDI for $\Lambda$-Fleming-Viot support}} {\textsc{Speed of CDI for $\Lambda$-Fleming-Viot support}}

\section{Introduction}
Coming down from infinity (CDI for short) is a fascinating phenomenon that manifests in various contexts for stochastic models. Intuitively, CDI occurs for a real-valued stochastic process when, starting with an infinite initial value at time $t = 0$, the process immediately assumes finite values at any positive time. In the context of real-valued Markov processes, CDI is linked to the concept that infinity serves as an entrance boundary.
Moreover, CDI is a significant aspect in the study of coalescent processes, which we will introduce in more detail in the subsequent section.

There is a growing  literature  on quantifying the speed at which CDI occurs in various stochastic models. For instance, Aldous \cite{MR1673235} determined the speed of CDI for the block counting process of the Kingman's coalescent with binary collisions.
Using an innovative martingale technique, Berestycki et al. \cite{BBL} demonstrated that when CDI occurs for the
$\Lambda$-coalescent with multiple collisions characterized by a coalescing mechanism denoted as $\psi$ (to be introduced in Section \ref{sec:2}) the block counting process $N(t)$ descends from infinity along a deterministic function $v$
that can be derived from $\psi$. Moreover, they established that the ratio $N(t)/v(t)$ converges both almost surely and in moments as $t \rightarrow0+$.
Additionally, Limic and Talarczyk \cite{MR3342667} explored the small-time second-order asymptotics for the block counting process $N(t)$ in this context. For more general $\Xi$-coalescents that permit simultaneous multiple collisions, Limic \cite{MR2594877} identified almost sure small-time asymptotics for $N(t)$ across a broad class of $\Xi$-coalescents that exhibit the CDI phenomenon.

Among the various results concerning CDI for different processes, Bansaye et al. \cite{MR3595771} studied the speed of CDI for birth-and-death processes by analyzing the short-time asymptotic behavior of the downcrossing time, rescaled by its expectation given the process starts at infinity. They established convergence in probability for gradual regime and convergence in distribution for fast regime.
In the context of diffusion processes that start from infinity and subsequently come down from infinity, Bansaye et al. \cite{MR4014294} proved further convergence results concerning hitting times after appropriate rescaling and renormalization. Their work demonstrated both a law of large numbers and a central limit theorem for these processes.
Additionally, Foucart et al. \cite{MR4255235} examined the CDI property for a class of nonlinear continuous-state branching processes defined as time-changed spectrally positive L\'evy processes stopped at $0$. They identified the speed of CDI in both slow and fast regimes.
These studies help to understand the dynamics involved in CDI across a variety of stochastic models.

The Fleming-Viot process is a probability measure-valued process that arises from population genetic models. It incorporates mechanisms such as reproduction (resampling) and mutation, where the measure is defined on a space of types. In this context, the reproduction is associated with a coalescent process, while the  mutation is often interpreted as spatial motion.
The probability measure in this context describes the relative frequencies of different types of individuals within the population. For comprehensive introductions to Fleming-Viot processes and related population models, we refer to Dawson \cite{Dawson93} and Etheridge \cite{Eth12}.
The Fleming-Viot process is dual to a coalescent, which captures the reproduction mechanism of the corresponding population model. In this paper, we focus on the $\Lambda$-Fleming-Viot process, which is dual to a $\Lambda$-coalescent and with Brownian spatial motion.

It is well established that a broad class of measure-valued processes, including the Fleming-Viot processes with Brownian spatial motion, exhibit what is known as the compact support property. This means that the measure-valued process has a compact support at any positive time, despite the initial measure having an unbounded support.  The earliest study on the compact support property of classical Fleming-Viot processes can be found in Dawson and Hochberg \cite{Dawson}. In contrast, if the spatial motion has jumps, the support can extend arbitrarily far away as time approaches zero; this is discussed in  Hughes and Zhou \cite{HZ23} for Fleming-Viot processes with L\'evy spatial motion.
The compact support property can also be viewed as a CDI-type phenomenon, raising a natural question about how to characterize the speed of this CDI. To the best of our knowledge, there has been no prior work addressing this.
The motivation of this paper is to further explore  the dynamics of the Fleming-Viot process with Brownian motion to understand the convergence behaviors associated with the CDI for its support.

As noted by Birkner et al. \cite{BB2} and Blath \cite{Jo}, if the dual
$\Lambda$-coalescent remains infinite, the support of the $\Lambda$-Fleming-Viot process will propagate to infinity instantaneously.
Therefore, we aim to study the CDI for the support process specifically by focusing on the $\Lambda$-Fleming-Viot process  for which the associated $\Lambda$-coalescent comes down from infinity.
We aim to characterize how the support  of the process evolves over time as it transitions from an initially unbounded set to a bounded set in different scenarios, and understand how the speed of CDI for  the
$\Lambda$-Fleming-Viot process is associated to the speed function $v$ for  $\Lambda$-coalescent.

The renowned lookdown representation for the Fleming-Viot process, initially introduced in the seminal wok of Donnelly and Kurtz \cite{DK96}, offers a particle system approximation to the measure-valued process. This particle system, along with its genealogical structure, is pivotal in investigation of the support properties of the Fleming-Viot process.
Exploiting the lookdown representation, Liu and Zhou \cite{LZ1} established a sufficient condition related to the rate of reduction for the block counting process, which is crucial for demonstrating the compact support property of the
$\Lambda$-Fleming-Viot process with Brownian spatial motion. Furthermore, they proved modulus of continuity results utilizing the lookdown representation in their subsequent works \cite{LZ2, LZ3}.

In this paper, we leverage the lookdown representation again to characterize the asymptotic behaviors of extreme values for the $\Lambda$-Fleming-Viot support.
To better illustrate our results, we focus on the one-dimensional case and assume that the upper bound of the initial support is infinite. Our goal is to quantify how quickly the upper bound of the $\Lambda$-Fleming-Viot support process decreases from infinity over a short time.

To characterize the speed of CDI for the support process using the speed function $v$ associated with the $\Lambda$-coalescent, we adapt the arguments presented in \cite{LZ2} to establish an integral test concerning $v$ for the modulus of continuity of the $\Lambda$-Fleming-Viot support process, which subsequently leads to the compact support property. In this result, we derive a relatively relaxed condition on $v$ compared to that in \cite{LZ2}, albeit at the expense of a more crude modulus of continuity. Nevertheless, this result is sufficient for us to specify the speed of the CDI for the support of the $\Lambda$-Fleming-Viot process with Brownian spatial motion. A moment estimate on the speed of CDI for the $\Lambda$-coalescent obtained in \cite{BBL} plays a crucial role in our proof.

The lookdown representation provides valuable insights on how the tail distribution $\bar{F}$ of the initial measure, the number of ancestors, and the dislocations of individuals at time $t$ from their respective ancestors at time $0$ contribute to this characterization.
Applying the modulus of continuity results for the Fleming-Viot support with Brownian spatial motion, we find that the dislocations of the descendants from their respective ancestors are sufficiently small near time $0$, rendering them negligible in affecting the speed of CDI. As a result, the speed of CDI is predominantly governed by the asymptotic behavior of the tail distribution of the initial probability measure, coupled with the speed of CDI for the corresponding $\Lambda$-coalescent, which underscores the interplay between  distribution of the initial measure and the genealogical structure determined by the coalescent.

When the spatial motion is first ignored, the lookdown representation further enables us to apply extreme value theory to find  the speed of the CDI for the ancestral process. To this end, we first establish a  convergence result for $M(t)$, the maximal position of the ancestors for the individuals at time $t$, rescaled using the speed function $v$ and the asymptotic inverse function of the tail distribution $\bar{F}$.
For Fleming-Viot process with Brownian spatial motion, by  modulus of continuity for  $\Lambda$-Fleming-Viot process we  show that
the aforementioned result still holds with $M(t)$ replaced by $\hat{M}(t)$ under a mild condition, where $\hat{M}(t)$, the extremal process, denotes the maximal position of the individuals alive at time $t$, i.e. the contribution of spatial motion to the speed of CDI is negligible.

More explicit convergence results can be derived for $\Lambda$-Fleming-Viot processes with  more specific asymptotic tail initial distribution $\bar{F}$.
If $\bar{F}$ behaves asymptotically like a power function (slow regime), we identify a function $a$ such that the rescaled process $a(v(t))\hat{M}(t)$ converges in distribution to a Fr\'echet distribution as $t \rightarrow0+$. Here, $a(v(t))^{-1}$ can be interpreted as the speed of the CDI. If $\bar{F}$ behaves like an exponential function asymptotically (fast regime), we find functions $\bar{a}$ and $\bar{b}$ such that the renormalized process $\bar{a}(v(t))\hat{M}(t) - \bar{b}(v(t))$ converges in distribution to the Gumbel distribution. Additionally, we show convergence in probability results for both regimes. Furthermore, these results can be generalized to $\Lambda$-Fleming-Viot processes on $\mathbb{R}^d $.

The remainder of the paper is organized as follows. Section \ref{sec:2} provides an introduction and preliminary results on the $\Lambda$-coalescent, including the speed of CDI for the $\Lambda$-coalescent, the $\Lambda$-Fleming-Viot process, and its lookdown representation. In this section, we also present an integral test on the modulus of continuity and discuss the compact support property of the $\Lambda$-Fleming-Viot process with Brownian spatial motion.
In Section \ref{sec:3}, we 
derive short-time convergence results in distribution and in probability, respectively, for the scaled or renormalized extremal processes.

Throughout this paper,  the letters $C$ or $c$, with or without subscripts,  denote constants whose exact values are not important and may vary from place to place. Given two functions $f, g: \mathbb{R}^{+} \rightarrow \mathbb{R}^{+}$,  write $f \sim g$ if
$\lim{f(x)}/{g(x)} = 1,$
and $f = o(g)$ if $\limsup{f(x)}/{g(x)} = 0,$ while $f = O(g)$ if
$\limsup{f(x)}/{g(x)} < \infty$,
where $x$ either goes to $0$ or $\infty$ in the limits
depending on the context. Denote by $f\circ g(\cdot)=f\(g\(\cdot\)\)$.

\section{Preliminaries}\label{sec:2}
\subsection{$\Lambda$-coalescent}
We refer to Bertoin \cite{Bertoin} for a systematic  introduction on exchangeable coalescents.
Let $[n]:=\{1,\ldots,n\}$ and $[\infty]:=\{1,2,\ldots\}$.
An ordered partition of a set $D \subset [\infty]$ is defined as a countable collection $\pi = \{ \pi_i, i=1, 2, \ldots \}$ of disjoint blocks such that the union of all blocks is equal to $D$, i.e., $\bigcup_{i} \pi_i = D$ and the minimum element of each block satisfies $\min \pi_i < \min \pi_j$ for $i < j $.
For convenience, we consider the case where $[\infty]$ is represented by the partition in which all blocks are merged into a single block. In this context, the blocks in $\pi$ are arranged in order based on their least elements.
We denote the set of ordered partitions of $[n]$ by $\mathcal{P}_n$, while $\mathcal{P}_{\infty}$ represents the set of ordered partitions of $[\infty]$.
Additionally, we define $\mathbf{0}_{[n]}$ as the partition of $[n]$ consisting of singletons, specifically $\mathbf{0}_{[n]} := \{ \{1\}, \{2\}, \ldots, \{n\} \} $. Similarly, $\mathbf{0}_{[\infty]}$ denotes the partition of $[\infty]$ that consists of singletons.

{\it Kingman's coalescent} is a time-homogeneous Markov process taking values in $\mathcal{P}_{\infty}$ where distinct pairs of blocks merge independently at the same rate. On the other hand, the $\Lambda$-coalescent extends Kingman's coalescent by allowing for multiple collisions between blocks. For detailed discussions, see Pitman \cite{Pitman99} and Sagitov \cite{Sag99}.
The $\Lambda$-coalescent is represented as a Markov process $\Pi := (\Pi(t))_{t \geq 0}$ that takes values in $\mathcal{P}_{\infty}$. For each $n \in [\infty]$, its restriction to $[n]$, denoted $\Pi_n := (\Pi_n(t))_{t \geq 0}$, behaves as a Markov process valued in $\mathcal{P}_n$. The transition rates for this process are defined as follows: given that there are $b$ blocks in the partition, each $k$-tuple of these blocks (where $2 \leq k \leq b$) merges independently to form a single block at a rate given by
$\lambda_{b,k}=\int_{[0,1]}x^{k-2}(1-x)^{b-k}\Lambda(\mathrm{d}x)$,
where $\Lambda$ is a finite measure on the interval $[0,1]$. For $n = 2, 3, \ldots$, we define the total coalescence rate starting with $n$ blocks as
$$\lambda_n=\sum_{k=2}^n{n\choose k}\lambda_{n,k}.$$
This expression sums over all possible $k$-tuples of $n$ blocks, providing a  rate for coalescence events among the blocks.

In particular, if we take $\Lambda(dx)=\Lambda\(\{0\}\)\delta_0(dx)$, the resulting coalescent process is Kingman's coalescent. On the other hand, if $\Lambda$ follows the Beta$(2-\be,\be)$ distribution, specifically defined as
\beqnn
\La(dx)=\frac{\Gamma(2)}{\Gamma(2-\be)\Gamma(\be)}x^{1-\be}\(1-x\)^{\be-1}dx
\eeqnn
for $\be\in(0,2)$, then the associated $\Lambda$-coalescent becomes a Beta$(2-\be,\be)$ coalescent.
If we consider $\Lambda(dx)=\Lambda\(\{1\}\)\delta_1(dx)$, the resulting coalescent process is known as the star-shaped coalescent. In this case, the transition rates satisfy  $\lambda_{b,k}=0$
for any $2\leq k\leq b-1$.
For the star-shaped coalescent, all blocks merge into a single block after an independent exponential time $T$ with a rate of $\Lambda(\{1\})$. This means that instead of allowing independent merging among pairs of blocks, the entire system collapses into one block at once after a random time dictated by the exponential distribution.
To prevent this scenario and ensure more complex coalescent behavior, we henceforth assume that $\Lambda$ has no atom at $1$. This condition allows for a richer structure in the coalescent process, avoiding the trivial case where all blocks merge instantaneously.

Given any $\La$-coalescent $\Pi=\(\Pi(t)\)_{t\geq 0}$ with $\Pi(0)=\mathbf{0}_{[\infty]}$, let $N(t)$ be the total number of blocks in the partition $\Pi(t)$. The $\La$-coalescent $\Pi$ {\it comes down from infinity} if $\mathbb{P}\(N(t)<\infty\)=1$ for all $t>0$ and it  {\it stays infinite} if $\mathbb{P}\(N(t)=\infty\)=1$ for all $t>0$.
If $\Lambda(\{1\})=0$, then the $\La$-coalescent either comes down from infinity with probability one or stays infinite with probability one; see Pitman \cite[Proposition 23]{Pitman99}.

\subsection{Speed function of CDI for $\Lambda$-coalescent}
Define
\beqnn
\psi(q)\ar:=\ar\int_{[0,1]}(e^{-qx}-1+qx)x^{-2}\Lambda(\mathrm{d}x)\cr
\ar=\ar\La(\{0\})\frac{q^2}{2}+\int_{(0,1]}(e^{-qx}-1+qx)x^{-2}\Lambda(\mathrm{d}x).
\eeqnn
Bertoin and Le Gall \cite{BerLeGa06} verified that a $\La$-coalescent comes down from infinity if and only if
\beqlb\label{eq:comin}
\int_a^\infty \frac{\mathrm{d}q}{\psi(q)}<\infty,
\eeqlb
where the integral is finite for some (and then for all) $a>0$.
Another necessary and sufficient condition for $\La$-coalescent to come down from infinity in terms of  rates of reduction for the block counting process can be found in Schweinsberg\cite{Jason2000}.

For the $\La$-coalescent coming down from infinity, let
$v: (0, \infty)\mapsto (0, \infty)$ be the continuous decreasing function satisfying
\begin{equation}\label{eq:vv}
\int_{v(t)}^\infty\frac{\mathrm{d}q}{\psi(q)}=t.
\end{equation}
For example, Kingman's coalescent with $\Lambda(dx)=\delta_0(dx)$ has a speed function $v(t)=2/t$ for all $t>0$. The
Beta$(2-\be,\be)$ coalescent with $\be\in(1,2)$ has a speed function $v(t)\sim\(\be\Gamma(\be)\)^{1/(\be-1)}t^{-1/(\be-1)}$ as $t\rightarrow0+$;
see \cite[Theorem 1.1]{BBS08AIHP}.

Berestycki et al. \cite{BBL} showed that
\beqlb\label{eq:cdi}
\lim_{t\rightarrow0+}\frac{N(t)}{v(t)}=1\text{\,\,\,\,a.s.}
\eeqlb
conditional on $\La([0,1))=1$.
The function $v$ is referred to as the speed of CDI for the corresponding $\Lambda$-coalescent.
A scaling of the total mass of $\Lambda$ by a constant factor  induces the scaling of the collision rates by the
same factor, and the speed $v$ of CDI  in (\ref{eq:vv}) above is also scaled in the same way.
For convenience 
we assume that $\Lambda([0,1])=1$.
Moreover, the following convergence in moments  was also shown in Berestycki et al. \cite{BBL}.
For any $\tilde{d}\in[1,\infty)$,
\beqlb\label{eq:cdimo}
\lim_{t\rightarrow0+}\mathbb{E}\(\sup_{s\in[0,t]}\Big|\frac{N(s)}{v(s)}-1\Big|^{\tilde{d}}\)=0.
\eeqlb

The next result is an adaptation of Berestycki et al. \cite[Proposition 12]{BBL}. Its proof is similar to that of \cite[Proposition 5.2]{LZ3}.
\begin{proposition}\label{pro:number}
If the associated $\Lambda$-coalescent $\Pi$ comes down from infinity, then for any $\al^{*}\in(0,1/2)$, we have
\beqnn
\mathbb{P}\({e^{-24  s^{\al^{*}}}{v(s)}\leq N\(s\)\leq e^{24  s^{\al^{*}}}{v(s)}}\)\geq 1-O\(s^{1-2\al^{*}}\)
\eeqnn
as $s\rightarrow0+$.
\end{proposition}

\subsection{$\Lambda$-Fleming-Viot process and lookdown representation}
In this subsection, we adopt the approach of Birkner and Blath \cite{BB} to present the $\Lambda$-Fleming-Viot process through its lookdown representation. The lookdown representation is a powerful technique that was first introduced in a series of papers by Donnelly and Kurtz \cite{DK96, DK99a, DK99b}. This representation provides a discrete particle framework for analyzing probability measure-valued Fleming-Viot processes, which are important in the study of Fleming-Viot processes.

For any collision measure $\Lambda$, let $\Lambda_0$ be its restriction on $(0,1]$. We proceed to define, on a filtered probability space $(\Omega, \mathcal{F}, \mathbb{P})$, $\(\mathbb{R}^d\)^{\infty}$-valued stochastic processes $\(X_1(t),X_2(t),\ldots\)_{t\geq 0}$ to describe the  dynamics of countably infinite many individuals in a population, and independent Poisson processes $\{\mathbf{N}_{ij}(t): 1\leq i<j<\infty\}$ with common arrival rate $\La(\{0\})$ together with an independent Poisson point process $\tilde{\mathbf{N}}$ on $\mathbb{R}_{+}\times\left(0, 1\right]$ with intensity measure $\mathrm{d}t\times x^{-2}\La_0\(\mathrm{d}x\)$ to describe the reproduction mechanism in the population.

In the particle system, $X_i(t)$ represents the spatial location (type) of the individual at ``level'' $i$ at time $t$ and the initial locations $\(X_1(0),X_2(0),\ldots\)$ are exchangeable random variables.
The particle system undergoes spatial movement for the individuals and  reproduction (resampling) that updates the genealogy of those individuals.

Single birth events occur at arrival times of $\{\mathbf{N}_{ij}(t): 1\leq i<j<\infty\}$. If $\Delta{\mathbf{N}}_{ij}(t)=1$, the individual attached level $i$ gives birth to a child, which starts at the same location as its mother and is attached level $j$. The original individuals with level $j$ or with levels higher than $j$ are shifted accordingly. More specifically,
\begin{eqnarray*}
X_k(t)\,=\,\begin{cases} X_k(t-), &\text{~~if $k<j$},\cr
X_i(t-), &\text{~~if $k=j$},\cr
X_{k-1}(t-), &\text{~~if $k>j$}.
\end{cases}
\end{eqnarray*}
Here, $i$ and $k-1$ ($k>j$) on the right hand side represent  ancestral levels for the individuals at level $j$ and level $k$, respectively, at time $t$.

Multiple birth events are determined by an independent Poisson point process $\tilde{\mathbf{N}}$ on $\mathbb{R}_{+}\times\left(0, 1\right]$.
At a jump point $\(t, x\)$, for each level we independently toss a coin with the probability of head being $x$ to determine the levels of newborns at time $t$.
Let $J_t^k$ be the collection of levels lower than or equal to $k$ with head appearing and let $J_t=J_t^\infty$. At this time, the individual at level $i=\min J_t$ gives birth to children at levels $J_t\setminus\{i\}$ at location of the mother.
The mother's level keeps unchanged. The original individuals attached levels in the set $[\infty]\setminus i$ are attached new levels from the set $[\infty]\setminus J_t$ but keep their original order, i.e., what used to be at a higher level is still at a higher level. More specifically, for any $k\in[\infty]$,
\begin{eqnarray*}
X_k(t)\,=\,
\begin{cases}
X_k(t-), \text{~for~} k\leq i,\\
X_i(t-), \text{~for~} k>i \text{~with~} k\in J_t,\\
X_{k-|J_t^k| +1}(t-), \text{~otherwise}.
\end{cases}
\end{eqnarray*}
Between jump times of the Poisson point process, individuals keep their levels unchanged and move independently according to standard Brownian motions in ${{\mathbb{R}}^d}$. Then $\(X_1(t),X_2(t),X_3(t),\ldots\)$ is exchangeable for each $t>0$.

Denote by $\mathcal{C}_b^2({{\mathbb{R}}^d})$ the collection of functions that are differentiable   up to the second order with bounded second derivative.
Define $M_1(\mathbb{R}^d)$ as the collection of probability measures on $\mathbb{R}^d$.
Let $f_1,\,f_2,\ldots$ be uniformly bounded functions in $\mathcal{C}_b^2({{\mathbb{R}}^d})$ that separate points in $M_1({{\mathbb{R}}^d})$ in the sense that $\int_{{{\mathbb{R}}^d}}f_k\mathrm{d}\mu=\int_{{{\mathbb{R}}^d}}f_k\mathrm{d}\nu$ for any $k\in[\infty]$ implies $\mu=\nu$, where $\mu,\,\nu\in M_1({{\mathbb{R}}^d})$; see e.g. \cite[Lemma 1.1]{DK96}.
Denote by $D([0, \infty), M_1({{\mathbb{R}}^d}))$ the space of c\`adl\`ag paths from $[0,\infty)$ to $M_1({{\mathbb{R}}^d})$.
Defined on $M_1({{\mathbb{R}}^d})$ a metric
\begin{equation*}
d\(\mu,\nu\):=\sum_k\frac{1}{2^k}\left|\int_{{{\mathbb{R}}^d}}f_k\mathrm{d}\mu-\int_{{{\mathbb{R}}^d}}f_k\mathrm{d}\nu\right|\text{~~for ~}\mu,\,\nu\in M_1({{\mathbb{R}}^d})
\end{equation*}
and equip $D([0, \infty),M_1({{\mathbb{R}}^d}))$ with the  topology of locally uniform convergence under metric
\beqnn
d_{p}\(\mu,\nu\):=\int_0^{\infty}e^{-t}d\(\mu(t),\nu(t)\)\mathrm{d}t.
\eeqnn

The $M_1({{\mathbb{R}}^d})$-valued process $\(X(t)\)_{t\geq 0}$, defined in terms of the particle system in the lookdown representation by
\beqnn
X(t):=\lim_{n\rightarrow\infty}X^{(n)}(t):= \lim_{n\rightarrow\infty}\frac{1}{n}\sum_{i=1}^n\delta_{X_i(t)},\,\,\,\,t\geq 0,
\eeqnn
is the $\La$-Fleming-Viot process with Brownian spatial motion. Moreover, the empirical process $\(X^{(n)}(t)\)_{t\geq0}$ converges almost surely on path space $D([0, \infty), M_1({{\mathbb{R}}^d}))$ to $\(X(t)\)_{t\geq0}$; see   \cite[Theorem 3.2]{DK99b}.  A similar result is shown in  Birkner et al. \cite[Theorem 1.1]{MJM} for $\Xi$-Fleming-Viot process with jump type spatial motion where the $\Xi$-coalescent generalizes the $\Lambda$-coalescent by allowing simultaneous multiple collisions.


Its generator, denoted by $\mathcal{L}$, acting on test functions of the form
$\mathbf{g}(\mu)=\langle\mu,g_1\rangle\langle\mu,g_2\rangle\cdots\langle\mu,g_n\rangle$ for $\mu\in M_1({{\mathbb{R}}^d})$, {$g_i\in \mathcal{C}_b^2\({{\mathbb{R}}^d}\)$,} $i\in[n]$, $n\in[\infty]$ and $\langle\mu,g_{i}\rangle =\int_{{{\mathbb{R}}^d}}g_i(x)\mu(\mathrm{d}x)$, is given by
\beqnn
\mathcal{L}\mathbf{g}(\mu)
\ar:=\ar \sum_{i=1}^{n}\langle\mu,\frac{\Delta}{2}g_i\rangle\prod_{\ell\neq i}\langle\mu,g_{\ell}\rangle\cr
\ar\ar
+\sum_{\substack{J\subseteq\{1,\ldots,n\}\\|J|\geq 2}}\lambda_{n,|J|}\(\langle\mu, \prod_{i\in J}g_i\rangle-\prod_{i\in J}\langle\mu,g_i\rangle\) \prod_{\ell\notin J}\langle\mu,g_{\ell}\rangle\cr
\ar=\ar \sum_{i=1}^{n}\langle\mu,\frac{\Delta}{2}g_i\rangle\prod_{\ell\neq i}\langle\mu,g_{\ell}\rangle \cr
\ar\ar +\Lambda(\{0\})\sum_{1\leq i<j\leq n}\(\langle\mu,g_i g_j\rangle-\langle \mu,g_i\rangle\langle\mu,g_j\rangle\)\prod_{\ell\neq i,\,j }\langle\mu,g_{\ell}\rangle\cr
\ar\ar+\int_{(0,1]} z^{-2}\Lambda\(dz\)\int_{\mathbb{R}^d} \[\mathbf{g}\(z\delta_{x}+\(1-z\)\mu\)-\mathbf{g}\(\mu\)\]\mu\(dx\),
\eeqnn
where {$\Delta$} is the Laplacian operator and $|J|$ is the cardinality of {$J$}.

The well-posedness of a martingale problem for $\Lambda$-Fleming-Viot process is shown in Donnelly and Kurtz \cite[Theorem 4.5]{DK99b}. Also see Li et al. \cite[Theorem 3.3]{llxz} for the well-posedness of a martingale problem for $\Xi$-Fleming-Viot process  with Brownian spatial motion.

Throughout the paper  we simply assume that $X_1(0),\,X_2(0),\ldots$ are i.i.d. random variables to facilitate the application of extreme values to the support process.  Additionally, we will assume that the $\Lambda$-coalescent comes down from infinity henceforth.

Write $S(\mu)$ for the closed support of measure $\mu$.
The following lemma
 is derived from \cite[Lemma 3.1]{LZ2} and \cite[Lemma 5.5]{LZ2}.
\begin{lemma}
For any $t\geq0$, we have $\mathbb{P}$-a.s.\,
$S\(X(t)\)=\overline{\{X_1(t),X_2(t),\ldots\}}.$
\end{lemma}

\subsection{$\Lambda$-coalescent in the lookdown representation}
For any $0\leq s\leq t$, denote by $\mathbb{L}_n^t(s)$ the left continuous ancestral level at time
$s$  for the individual with level $n$ at time $t$.  We refer to \cite[equation (1.45)]{BB} for the specific expression of $\mathbb{L}_n^t(s)$.
Write $\big(\Pi^t(s)\big)_{0\leq s\leq t}$ for the
$\mathcal{P}_\infty$-valued process such that $i$ and $j$ belong to
the same block of $\Pi^t(s)$ if and only if $\mathbb{L}_i^t(t-s)=\mathbb{L}_j^t(t-s)$.
The process $\big(\Pi^t(s)\big)_{
0\leq s\leq t}$ has the same law as the
$\La$-coalescent running up to time $t$; see \cite[Section 5]{DK99b}  or \cite[Section 1.3]{BB}.
For any $0\leq s\leq t$, write
\beqnn
N({s,t})\ar:=\ar\#\Pi^{t}\(t-s\)
\eeqnn
for the number of blocks in $\Pi^{t}\(t-s\)$ and
$\Pi^{t}\(t-s\)=\{\pi^{s,t}_{\ell}: 1\leq \ell\leq N({s,t})\},$
where $\pi^{s,t}_{\ell}, 1\leq \ell\leq N({s,t})$ are all the
disjoint blocks ordered by their least
elements. Since $N(0,t)$ will be used repeatedly in the following, for simplicity, we abuse notation by writing $N(t)$ instead of $N(0,t)$ below. They are not identical but have the same distribution.

A result related to the ancestral level and the coalescent process is stated below, as given in \cite[Lemma 3.1]{LZ1}.
\begin{lemma}\label{level}
For any fixed $T>0$, let $\(\Pi^{T}(t)\)_{0\leq t\leq T}$ be
the $\La$-coalescent recovered from the lookdown
representation. Then given any $t\in[0,T]$ and the ordered random
partition $\Pi^T(t)=\left\{\pi_{\ell}(t):
\ell=1,\ldots,\#\Pi^T(t)\right\}$, we have \,
$\mathbb{L}_j^T\(T-t\)=\ell\text{~~for any~~}j\in\pi_{\ell}(t).$
\end{lemma}

\subsection{Ancestral process }
For any $T>0$,
denote by
$\(X_{1,t},X_{2,t},X_{3,t},\ldots\)_{0\leq t\leq T}$
the {\it ancestral process} recovered from the lookdown representation by
\beqnn
X_{i,t}\(s\):=
X_{\mathbb{L}_i^t\(s-\)}{\(s-\)} \text{~~for~~}0\leq s\leq t.
\eeqnn
By the definition, $\(X_{i,t}(s)\)_{0\leq s\leq t}$ keeps a record of the trajectory for the $i$-th individual alive at time $t$
and $\(X_{1,t},X_{2,t},X_{3,t},\ldots\)$ depicts the genealogy of all those individuals alive at time $t$.
We call $\(X_{1,t},X_{2,t},X_{3,t},\ldots\)$ the ancestral process recovered at time $t$.

For any $0\leq r< s\leq t$, let $H^{t}(r,s)$ be the maximal dislocation  between times $r$ and $s$  of the ancestral process
recovered at time $t$, i.e.
\beqnn
H^{t}(r,s)\ar:=\ar\max_{j\in [\infty]}\left\lvert {X_{\mathbb{L}_j^{t}(s-)}{(s-)}-X_{\mathbb{L}_j^{t}(r-)}{\(r-\)}}\right\rvert\cr
\ar=\ar\max_{1\leq \ell\leq
N(r,t)}\max_{j\in\pi^{r,t}_{\ell}}\left\lvert {X_{\mathbb{L}_j^{t}(s-)}{\(s-\)}-X_{\ell}{(r-)}}\right\rvert,
\eeqnn
where Lemma \ref{level} is applied in the last line.

\subsection{Modulus of continuity for $\Lambda$-Fleming-Viot process with Brownian spatial motion}
Previous results on the modulus of continuity for the $\Lambda$-Fleming-Viot process with Brownian spatial motion were established in \cite{LZ2, LZ3} under specific conditions regarding the coalescing rates of the $\Lambda$-coalescent. In this paper, we aim to express the speed of CDI in terms of $v$; therefore, it would be more convenient to have such a result under conditions related to $v$. We thus propose the following condition.

\begin{assumption}\label{weaker_con}
There  exist some constants $\delta_o\in(0,1)$ and $\ep_o>0$ such that
\beqlb\label{eq:v1}
\int_{0+}v(t)^{1+\ep_o}e^{-t^{-\delta_o}}\mathrm{d}t<\infty.
\eeqlb
\end{assumption}
Note that any power speed  function $v$ satisfies the above condition. For example, the $v$ for Beta$(\be,2-\be)$-coalescent with $\be\in(1,2]$ satisfies (\ref{eq:v1}).
We next present a modulus of continuity result, with relaxed condition compared to that in \cite{LZ2},  on the speed of CDI for the associate $\Lambda$-coalescent below.  Then the uniform compactness for the support follows as a corollary.
The proofs are similar to those in \cite{LZ2} and are deferred to the Appendix.

\begin{proposition}\label{le:weakm}
Under Assumption \ref{weaker_con},
for any $\nu\in M_1({{\mathbb{R}}^d})$, $T>0$ and $\delta\in(0,{(1-\delta_o)/2})$,  there exist a positive random variable $\Delta$ and a constant $C(\delta)$ such that $\mathbb{P}_{\nu}$-a.s. for all $0\leq s<t\leq T$ satisfying $t-s\leq\Delta$,
$$H^t\({s,t}\)\leq C(\delta)(t-s)^{\delta}.$$
Moreover,
 for any $x>0$ and $d_o>1/\ep_o$,
\beqlb\label{eq:Delta}
\mathbb{P}\(\Delta\leq x\)
\ar\leq\ar
\int_{0}^{2^{-\lceil\log(1/x)/\log 2\rceil}}v(t)^{1+\ep_o}e^{-t^{-\delta_o}}\mathrm{d}t+\frac{2^{-\lceil\log(1/x)/\log 2\rceil\(d_{o}\ep_{o}-1\)}}{1-2^{-(d_o\ep_o-1)}}:= h(x),
\eeqlb
where $\lceil y\rceil:=\min\{k\geq y, k\in[\infty]\}$.
\end{proposition}

\begin{corollary}\label{co:com}
Under Assumption \ref{weaker_con},
$\mathbb{P}$-a.s. the support  $S\(X(t)\)$ is  compact for all $t>0$.
\end{corollary}

\begin{remark}
To compare the integral tests (\ref{eq:comin}) and (\ref{eq:v1}), note that if $\psi(q)\sim Cq\(\log q\)^{1+\ep}$ as $q\rightarrow\infty$ for some $\ep>1$, then  $\psi(q)$ satisfies both (\ref{eq:comin}) and (\ref{eq:v1}). However, if
$\ep\in(0,1]$, then it satisfies (\ref{eq:comin}) but not (\ref{eq:v1}).
Sufficient conditions in terms of coalescing rates  are also found in \cite{LZ1} and \cite{LZ2} for the compact support property of $\Lambda$-Fleming-Viot processes with Brownian spatial motion. In particular, the compact support property holds if
$\La([0, x])\geq Cx^{\gamma}$ in a right neighborhood of $0$ for $\ga\in(0,1)$.
By a Tauberian argument, if $\La([0, x])\sim Cx^{\gamma}$ as $x\rightarrow0+$, then $\psi(q)\sim Cq^{2-\ga}$ as $q\rightarrow\infty$.
Such conditions thus imply (\ref{eq:v1}).
\end{remark}

\section{Speed of CDI for $\Lambda$-Fleming-Viot process}\label{sec:3}
The speed of CDI for the $\Lambda$-Fleming-Viot process is meaningful
only under the condition that the associated $\Lambda$-coalescent  comes down from infinity.
Therefore, throughout Section \ref{sec:3} we always assume that the associated $\Lambda$-coalescent comes down from infinity
with the speed of CDI $v(t)$ given by (\ref{eq:vv}).
Additionally, let $X_1(0),\,X_2(0),\ldots$ be i.i.d. random variables with  a common distribution function $F$ on $\mathbb{R}$.
We denote the probability law of the associated 
Fleming-Viot process by $\mathbb{P}_F$. In the following discussions, we will simply use $\mathbb{P}$ when there is no risk of confusion.

\subsection{Convergence in distribution for the maximal location of ancestors}
For all the individuals alive at time $t$ in the lookdown representation, define $(M(t))_{t\geq 0}$ as the maximal location of their ancestors at time $0$. By Lemma \ref{level},  $M(t):=\max \{X_i(0),1\leq i\leq N(t)\}, \, t>0$, which decreases in $t$.

 To characterize the asymptotics of $M(t)$ as $t\to 0+$, we next introduce an asymptotic inverse function for the tail distribution $\bar{F}:=1-F$; see Bingham et al. \cite[Theorem 1.5.12]{MR1015093} for a similar definition. Not that the asymptotic inverse function $\bar{F}^{(-1)}$ is not unique.

\begin{definition}\label{de:1}
We call a nonincreasing function $\bar{F}^{(-1)}$  an asymptotic inverse of nonincreasing function $\bar{F}$ if
\beqnn
\bar{F}\circ\bar{F}^{(-1)}(y)= \bar{F}\(\bar{F}^{(-1)}(y)\)\sim y\,\,\,\,\text{as\,\,\,}y\rightarrow0+,
\eeqnn
where $\sim$ denotes asymptotic equivalence. This means that $\bar{F}^{(-1)}(y)$ describes the behavior of the inverse of $\bar{F}(y)$
in the limit as $y$ approaches $0$ from above.
\end{definition}

We first present a scaling limit result on $M$ that does not depend on the spatial motion.

\begin{lemma}\label{th:uc}
For any $x>0$, we have
\beqnn
\ar\ar\lim_{t\rightarrow0+}
\mathbb{P}\Big(M(t)\leq\bar{F}^{(-1)}\(x/v(t)\)\Big)=e^{-x}.
\eeqnn
\end{lemma}

\begin{proof}
The condition $\bar{F}\circ\bar{F}^{(-1)}(y)\sim y$ as $y\rightarrow0+$
implies that given any $\ep>0$, we have
\beqnn
\({1-\ep}\)x/v(t)\leq \bar{F}\circ\bar{F}^{(-1)}(x/v(t))\leq \({1+\ep}\)x/v(t)
\eeqnn
for any $x>0$ as $t\rightarrow0+$, which yields
\beqnn
1-\({1+\ep}\)x/v(t)\leq {F}\circ\bar{F}^{(-1)}(x/v(t))\leq 1-\({1-\ep}\)x/v(t).
\eeqnn
Then we obtain
\beqnn
\ar\ar\mathbb{P}\big(M(t)\leq \bar{F}^{(-1)}\(x/v(t)\)\big)=\mathbb{E}\Big[\big(F\circ\bar{F}^{(-1)}\(x/v(t)\)\big)^{N(t)}\Big]
\leq \mathbb{E}\[\(1-(1-\ep)x/v(t)\)^{N(t)}\]
\eeqnn
and
\beqnn
\ar\ar\mathbb{P}\(M(t)\leq \bar{F}^{(-1)}\(x/v(t)\)\)\geq \mathbb{E}\[\(1-(1+\ep)x/v(t)\)^{N(t)}\].
\eeqnn
 It follows from Proposition \ref{pro:number} that as $t\rightarrow0+$,
\beqnn
\mathbb{P}\(N(t)< (1-\ep)v(t)\)\leq O(t^{1-2\al^{*}})
\eeqnn
and
\beqnn
\mathbb{P}\(N(t)\leq (1+\ep)v(t)\)\geq 1-O(t^{1-2\al^{*}}).
\eeqnn
Therefore, we continue to derive
\beqnn
\ar\ar
\limsup_{t\rightarrow 0+}\mathbb{P}\(M(t)\leq \bar{F}^{(-1)}\(x/v(t)\)\)\cr
\ar\ar\leq\limsup_{t\rightarrow 0+}\mathbb{E}\[\(1-(1-\ep)x/v(t)\)^{N(t)}1_{\{N(t)\geq (1-\ep)v(t)\}}\]+\mathbb{P}\(N(t)< (1-\ep)v(t)\)\cr
\ar\ar\leq\lim_{t\rightarrow 0+}\(1-\({1-\ep}\)x/v(t)\)^{(1-\ep)v(t)}+O(t^{1-2\al^{*}})\cr
\ar\ar=e^{-\(1-\ep\)^2x}
\eeqnn
and
\beqnn
\ar\ar\liminf_{t\rightarrow 0+}\mathbb{P}\(M(t)\leq \bar{F}^{(-1)}\(x/v(t)\)\)\cr
\ar\geq\ar\liminf_{t\rightarrow 0+}{\mathbb{E}\[\(1-(1+\ep)x/v(t)\)^{N(t)}1_{\{N(t)\leq (1+\ep)v(t)\}}\]}\cr
\ar\geq\ar\lim_{t\rightarrow 0+}\(1-\({1+\ep}\)x/v(t)\)^{(1+\ep)v(t)}\times\(1-O(t^{1-2\al^{*}})\)\cr
\ar=\ar e^{-\(1+\ep\)^2x}.
\eeqnn
Letting $\ep\rightarrow0+$, we obtain the desired result.
\end{proof}

\subsection{Convergence in distribution for the extremal process}
Let $\hat{M}(t)$, called the extremal process, be the maximal location of all the individuals alive at $t$, i.e. ${\hat{M}}(t):=\sup\{X_i(t),i=1,2,\ldots\}$.
Write $\hat{S}(t):=\sup S(X(t))$ for the  maximum of the closed support for $X(t)$.	Recall that for any $t\geq 0$, $\mathbb{P}$-a.s.
$S(X(t))= \overline{\{X_1(t),X_2(t),\ldots\}}$, which implies $\hat{M}(t)=\hat{S}(t) \,\,\mathbb{P}$-a.s. for any $t>0$.
The following corollary is an immediate consequence of Proposition \ref{le:weakm}.
\begin{corollary}\label{co:mo}
Under Assumption \ref{weaker_con} and for the random variable $\Delta$ and the constants $\delta$ and $C(\delta)$ in  Proposition \ref{le:weakm}, we have $\mathbb{P}$-a.s.
\beqlb\label{eq:mod}
M(t)-{C(\delta)t^{\delta}}\leq {\hat{M}}(t)\leq M(t)+{C(\delta)t^{\delta}}, \, 0<t\leq\Delta.
\eeqlb
Moreover,
\beqlb\label{eq:h}
\mathbb{P}\(\Delta> t\)=1-\mathbb{P}\(\Delta\leq t\)\geq 1-h(t)
\eeqlb
where function $h$ is defined in (\ref{eq:Delta}) satisfying $\lim_{t\rightarrow0+}h(t)=0$.
\end{corollary}
Equation (\ref{eq:mod}) provides an upper bound for the displacement of individuals relative to their respective ancestors over a sufficiently small time $t$. This bound plays a crucial role in demonstrating that the spatial motion does not affect the speed of the CDI.  We assume throughout the rest of the paper that Assumption \ref{weaker_con} holds.

\begin{theorem}\label{th:uc2}
Given the constants $\delta$ and $C(\delta)$ in \eqref{eq:mod} and  for any $x>0$, if
\beqlb\label{eq:Fa}
\lim_{t\rightarrow0+}v(t)\big|\bar{F}\(\bar{F}^{(-1)}\(x/v(t)\)\pm{C(\delta)t^{\delta}}\)-x/v(t)\big|=0,
\eeqlb
we have
\beqnn
\lim_{t\rightarrow0+}
\mathbb{P}\Big(\hat{M}(t)\leq\bar{F}^{(-1)}\(x/v(t)\)\Big)=e^{-x}.
\eeqnn
\end{theorem}
\begin{proof}
Given any $\ep\in(0,1)$, condition (\ref{eq:Fa}) implies that as $t\rightarrow0+$,
\beqnn
(x-\ep) /v(t)\leq\bar{F}\(\bar{F}^{(-1)}\(x/v(t)\)\pm{C(\delta)t^{\delta}}\)\leq (x+\ep)/v(t)
\eeqnn
for any $x>0$. Consequently,
\beqnn
1-\(x+\ep\)/v(t)\leq{F}\(\bar{F}^{(-1)}\(x/v(t)\)\pm{C(\delta)t^{\delta}}\)\leq1-\(x-\ep\)/v(t).
\eeqnn
For fixed $t>0$, we have
\beqnn
\mathbb{P}\(\hat{M}(t)\leq \bar{F}^{(-1)}\(x/v(t)\)\)\leq \mathbb{P}\(\hat{M}(t)\leq \bar{F}^{(-1)}\(x/v(t)\), t< \Delta\)+\mathbb{P}\(\Delta\leq t\).
\eeqnn
Applying (\ref{eq:mod}) and (\ref{eq:h}), we can derive
\beqnn
\mathbb{P}\(\hat{M}(t)\leq \bar{F}^{(-1)}\(x/v(t)\)\)
\ar\leq\ar\mathbb{P}\(M(t)\leq \bar{F}^{(-1)}\(x/v(t)\)+{C(\delta)t^{\delta}}\)+h(t)\\
\ar=\ar\mathbb{E}\[\({F}\big(\bar{F}^{(-1)}\(x/v(t)\)+{C(\delta)t^{\delta}}\big)\)^{N(t)}\]+h(t)\\
\ar\leq\ar\mathbb{E}\[\(1-(x-\ep)/v(t)\)^{N(t)}\]+h(t)
\eeqnn
and
\beqnn
\mathbb{P}\(\hat{M}(t)\leq \bar{F}^{(-1)}\(x/v(t)\)\)\ar\geq\ar\mathbb{P}\(\hat{M}(t)\leq \bar{F}^{(-1)}\(x/v(t)\), t< \Delta\)\\
\ar\geq\ar\mathbb{E}\[\(1-\(x+\ep\)/v(t)\)^{N(t)}\]\times \mathbb{P}\(t< \Delta\)\\
\ar\geq\ar\mathbb{E}\[\(1-\(x+\ep\)/v(t)\)^{N(t)}\]\times\(1-h(t)\)
\eeqnn
for any $x>0$.
Note that $h(t)$ tends to $0$ as $t\rightarrow0+$, thus, $h(t)$ has no effect on estimating the limit as
$t$ approaches $0$, so it can be ignored.
Similar to the proof of Lemma \ref{th:uc}, the conclusion will be proven based on the fluctuation of $N(t)$ around $v(t)$ by employing the law of total expectation to refine and bound the result. Since the proofs of the two are very similar, we omit the subsequent steps.
\end{proof}

\begin{remark}
	Observe that if
	$$
	\mathbb{P}\Big(\hat{M}(t)\leq\bar{F}^{(-1)}\(x/v(t)\)\Big)\to e^{-x}
	$$ as $t\to 0+$,
	then
	$$\mathbb{P}\Big(\bar{F}(\hat{M}(t))\geq\bar{F}(\bar{F}^{(-1)}\(x/v(t)))=x/v(t)-o(1/v(t)\)\Big)\to e^{-x}, $$
	which implies
	$$\mathbb{P}\Big(v(t)\bar{F}(\hat{M}(t))\geq x \Big)\to e^{-x}$$
 as $t\to 0+$.
\end{remark}

We first consider a case that $F$ is very heavy tailed that leads to very slow speed of CDI for the initial support.
\begin{corollary}
If $\bar{F}(x)\sim r_1(\log x)^{-r_3}$ as $x\rightarrow\infty$ for $r_1, r_3>0$, then for any $x\in\mathbb{R}$,
$$\lim_{t\rightarrow0+}\mathbb{P}\Big(r_3\log\log \hat{M}(t)-\log v(t)-\log r_1\leq x \Big)=e^{-e^{-x}},$$
i.e. the limiting distribution is the Gumbel distribution.
\end{corollary}

\begin{proof}

Let $g(y)=e^{\({y}/{r_1}\)^{-1/r_3}}$ as $y\rightarrow0+$. Given any $\ep>0$, it is straightforward to verify that
	\beqnn
	\bar{F}\(g(y)\)\leq(1+\ep){r_1}\(\log g(y)\)^{-r_3}=(1+\ep)y
	\eeqnn
	and
	\beqnn
	\bar{F}\(g(y)\)\geq(1-\ep){r_1}\(\log g(y)\)^{-r_3}=(1-\ep)y
	\eeqnn
	as $y\rightarrow0+$. Thus $\bar{F}\(g(y)\)\sim y$ as $y\rightarrow0+$. This implies that $\bar{F}^{(-1)}(y)=g(y)$ is an asymptotic inverse function for $\bar{F}(x)$.

In addition, for any constants $\delta, C>0$ we have
\beqnn
\ar\ar\bar{F}\(\bar{F}^{(-1)}\(x/v(t)\)\pm{Ct^{\delta}}\)\\
\ar\ar\leq(1+\ep)r_1\[\log\(\bar{F}^{(-1)}\(x/v(t)\)\pm{Ct^{\delta}}\)\]^{-r_3}\\
\ar\ar\leq(1+\ep)r_1\[\log\(\bar{F}^{(-1)}\(x/v(t)\)\)+\log\(1\pm{Ct^{\delta}}/\bar{F}^{(-1)}\(x/v(t)\)\)\]^{-r_3}\\
\ar\ar\leq(1+\ep)^2r_1\[\log\(\bar{F}^{(-1)}\(x/v(t)\)\)\]^{-r_3}\\
\ar\ar=(1+\ep)^2x/v(t)
\eeqnn
as $t\rightarrow0+$ and similarly
\beqnn
\ar\ar\bar{F}\(\bar{F}^{(-1)}\(x/v(t)\)\pm{Ct^{\delta}}\)\geq(1-\ep)^2x/v(t)
\eeqnn
as $t\rightarrow0+$. Therefore, condition (\ref{eq:Fa}) is satisfied.
By Theorem \ref{th:uc2}, for any $x>0$,
\beqnn
e^{-x}\ar=\ar\lim_{t\rightarrow0+}\mathbb{P}\Big(\hat{M}(t)\leq e^{\(x/(v(t){r_1})\)^{-1/r_3}}\Big)\\
\ar=\ar\lim_{t\rightarrow0+}\mathbb{P}\Big(r_3\log\log \hat{M}(t)-\log v(t)-\log r_1\leq -\log x \Big).
\eeqnn
Then the desired result follows by a change of variable.
\end{proof}

\subsection{Exact speed of CDI in slow regime}\label{slow_regime}
To obtain more explicit scaling limits for the  extremal process, we  impose the following asymptotic condition.
\beqlb\label{eq:pluplow}
\lim_{x\rightarrow\infty}\frac{\bar{F}(x)}{x^{-r_2}\(\log x\)^{r_3}}= r_1,
\eeqlb
where  $r_1>0,\,r_2>0$, $r_3\in\mathbb{R}$ are constants.
Further define function
$$a(x)\equiv a(x;r_1,r_2,r_3):=r_1^{-r_2^{-1}}r_2^{r_3r_2^{-1}}x^{-r_2^{-1}}\(\log x\)^{-r_3r_2^{-1}}
\quad\text{for}\quad x>1.$$

\begin{theorem}\label{le:e1}
Under condition (\ref{eq:pluplow}), for any $x>0$ we have
$$\lim_{t\rightarrow0+} \mathbb{P}\(a(v(t)) \hat{M}(t)\leq x\) =e^{-x^{-{r_2}}}.$$
\end{theorem}

\begin{proof}
Under condition (\ref{eq:pluplow}), one can show that
$$\bar{F}^{(-1)}(y)=r_1^{r_2^{-1}}r_2^{-r_3r_2^{-1}}y^{-r_2^{-1}}\(\log y^{-1}\)^{r_3r_2^{-1}} $$
is an asymptotic inverse for $\bar{F}$ as $y\rightarrow0+$.

Given any $x>0$, as $t\to 0+$, we have
\beqnn
\ar\ar\bar{F}\(\bar{F}^{(-1)}\(x/v(t)\){\pm}{C(\delta)t^{\delta}}\)\\
\ar\ar\sim r_1\(\bar{F}^{(-1)}\(x/v(t)\){\pm}{C(\delta)t^{\delta}}\)^{-r_2}\(\log\(\bar{F}^{(-1)}\(x/v(t)\){\pm}{C(\delta)t^{\delta}}\)\)^{r_3}\\
\ar\ar \sim r_1\[\(\bar{F}^{(-1)}\(x/v(t)\)\)^{-r_2}{\mp}r_2{C(\delta)t^{\delta}}\(\bar{F}^{(-1)}\(x/v(t)\)\)^{-r_2-1}\]\(\log\(\bar{F}^{(-1)}\(x/v(t)\)\)\)^{r_3}\\
\ar\ar\sim x/v(t){\mp}r_1r_2{C(\delta)t^{\delta}}\(\bar{F}^{(-1)}\(x/v(t)\)\)^{-r_2-1}\(\log\(\bar{F}^{(-1)}\(x/v(t)\)\)\)^{r_3}\\
\ar\ar\sim x/v(t)\mp r_1^{-r_2^{-1}}r_2^{1+r_3r_2^{-1}}{C(\delta)t^{\delta}}\(v(t)/x\)^{-1-r_2^{-1}}\(\log v(t)\)^{-r_3r_2^{-1}},
\eeqnn
which yields
\beqlb\label{eq:barF}
\lim_{t\rightarrow0+}v(t)\left|{\bar{F}\(\bar{F}^{(-1)}\(x/v(t)\){\pm}{C(\delta)t^{\delta}}\)-x/v(t)}\right|=0.
\eeqlb
Thus, condition (\ref{eq:Fa}) holds.

By Theorem \ref{th:uc2}, as $t\rightarrow0+$, we have
\beqlb\label{eq:u1}
\lim_{t\rightarrow0+}\mathbb{P}\Big(\hat{M}(t)\leq \bar{F}^{(-1)}(x/v(t))\Big)=e^{-x}
\eeqlb
for any $x>0$.

Given any $\ep>0$, as $t\to 0+$, it is straightforward to verify that
\beqnn
\bar{F}^{(-1)}\((x+\ep)/v(t)\)
\ar\leq \ar r_1^{r_2^{-1}}r_2^{-r_3r_2^{-1}}(x/v(t))^{-r_2^{-1}}\(\log v(t)\)^{r_3r_2^{-1}}
\leq\bar{F}^{(-1)}\((x-\ep)/v(t)\).
\eeqnn
It follows from (\ref{eq:u1}) that
\beqnn
e^{-(x+\ep)}\ar=\ar\lim_{t\rightarrow0+}\mathbb{P}\Big(\hat{M}(t)\leq \bar{F}^{(-1)}((x+\ep)/v(t))\Big)\\
\ar\leq\ar\liminf_{t\rightarrow0+}\mathbb{P}\Big(\hat{M}(t)\leq r_1^{r_2^{-1}}r_2^{-r_3r_2^{-1}}(x/v(t))^{-r_2^{-1}}\(\log v(t)\)^{r_3r_2^{-1}}\Big)\\
\ar\leq\ar\limsup_{t\rightarrow0+}\mathbb{P}\Big(\hat{M}(t)\leq r_1^{r_2^{-1}}r_2^{-r_3r_2^{-1}}(x/v(t))^{-r_2^{-1}}\(\log v(t)\)^{r_3r_2^{-1}}\Big)\\
\ar\leq\ar\lim_{t\rightarrow0+}\mathbb{P}\Big(\hat{M}(t)\leq \bar{F}^{(-1)}((x-\ep)/v(t))\Big)=e^{-(x-\ep)}.
\eeqnn
Letting $\ep\to 0+$ we have
\beqnn
\lim_{t\rightarrow0+}\mathbb{P}\Big(\hat{M}(t)\leq r_1^{r_2^{-1}}r_2^{-r_3r_2^{-1}}(x/v(t))^{-r_2^{-1}}\(\log v(t)\)^{r_3r_2^{-1}}\Big)= e^{-x}.
\eeqnn
The desired result follows by a change of variable.
\end{proof}

\begin{remark}\label{re:1}
The key to the proof for Theorem \ref{le:e1} is to identify the asymptotic inverse function for $\bar{F}$ that can be carried out through the following heuristic argument.
Assuming that for $x\to\infty$, function $y$ of $x$ satisfies
\beqlb\label{eq:yslow}
y\sim r_1x^{-r_2}\(\log x\)^{r_3},
\eeqlb
then
\beqnn
\log y\sim\log r_1-r_2\log x+r_3\log \log x.
\eeqnn
 Consequently, we anticipate
\beqlb\label{a_inv}
\log x\sim r_2^{-1}\log \(y^{-1}\)+\log C_1+C_2\log \log \(y^{-1}\)
\eeqlb
as $y\to 0+$. Substituting the above back into (\ref{eq:yslow}) allows to identify the constants
\beqnn
C_1=r_1^{r_2^{-1}}r_2^{-r_3r_2^{-1}}\,\,\,\,\text{and}\,\,\,\, C_2=r_3r_2^{-1}.
\eeqnn
Therefore, solving (\ref{a_inv}) for $x$ we expect
\beqnn
\bar{F}^{(-1)}(y)\sim r_1^{r_2^{-1}}r_2^{-r_3r_2^{-1}}y^{-r_2^{-1}}\(\log y^{-1}\)^{r_3r_2^{-1}}
\eeqnn
as $y\to 0+$.
\end{remark}

\begin{corollary}\label{co:re}
	If 
	$\bar{F}(x)\sim {r_1}x^{-r_2}$ as $x\rightarrow\infty$ for $r_1,\,r_2>0$, then
	$r_1\hat{M}(t)^{-r_2}v(t)\to \mathtt{e}$ in distribution as $t\rightarrow0+$, where $\mathtt{e}$ follows an exponential distribution with parameter $1$.
\end{corollary}
\begin{proof}
The asymptotic tail distribution satisfies condition (\ref{eq:pluplow}) by setting $r_3=0$.
Applying Theorem \ref{le:e1}, for any $x>0$, we derive
\beqnn
e^{-x^{-{r_2}}}\ar=\ar\lim_{t\rightarrow0+} \mathbb{P}\(r_1^{-r_2^{-1}}v(t)^{-r_2^{-1}}\hat{M}(t)\leq x\)=\lim_{t\rightarrow0+} \mathbb{P}\(r_1\hat{M}(t)^{-r_2}v(t)\geq x^{-r_2}\).
\eeqnn
Then the desired result follows by a change of variable.
\end{proof}

If the equality in condition (\ref{eq:pluplow}) is replaced by inequalities, then we can similarly show  the upper and lower bounds, respectively, for the scaling limits of $\hat{M}$.	

\begin{proposition}\label{pro:ere}
	\begin{enumerate}[(a)]
\item\label{a21}
If
\beqlb\label{eq:plup}
\limsup_{x\rightarrow\infty}\frac{\bar{F}(x)}{x^{-r_2}\(\log x\)^{r_3}}\leq r_1,
\eeqlb
then for any $x>0$, we have
\beqnn
\liminf_{t\rightarrow0+}\mathbb{P}\(a(v(t))\hat{M}(t)\leq x\){\geq}e^{-x^{-{r_2}}}.
\eeqnn
\item\label{a22}
If
\beqlb\label{eq:pllow}
\liminf_{x\rightarrow\infty}\frac{\bar{F}(x)}{x^{-r_2}\(\log x\)^{r_3}}\geq r_1,
\eeqlb
then for any $x>0$, we have
\beqnn
\limsup_{t\rightarrow0+}\mathbb{P}\(a(v(t))\hat{M}(t)\leq x\){\leq}e^{-x^{-{r_2}}}.
\eeqnn
\end{enumerate}
\end{proposition}

\begin{proof}
(a) Let $\bar{G}(x)=r_1x^{-r_2}\(\log x\)^{r_3}$ as $x\rightarrow\infty$. It is straightforward to verify that as $y\rightarrow0+$
$$\bar{G}^{(-1)}(y)=r_1^{r_2^{-1}}r_2^{-r_3r_2^{-1}}y^{-r_2^{-1}}\(\log \(y^{-1}\)\)^{r_3r_2^{-1}}$$
is an asymptotic inverse function for $G$. Given any $\ep>0$, as $t\rightarrow0+$, we have
\beqnn
\bar{G}^{(-1)}\((x+\ep)/v(t)\)\leq x^{-r_2^{-1}}/a(v(t))
\eeqnn
for any $x>0$.
This fact, combined with Corollary \ref{co:mo}, yields
\beqnn
\mathbb{P}\(a(v(t))\hat{M}(t)\leq x^{-r_2^{-1}}\)\ar\geq\ar \mathbb{P}\(\hat{M}(t)\leq \bar{G}^{(-1)}\((x+\ep)/v(t)\)\)\\
\ar\geq\ar\mathbb{P}\(\hat{M}(t)\leq \bar{G}^{(-1)}((x+\ep)/v(t)),t<\Delta\)\\
\ar\geq\ar\mathbb{P}\({M}(t)\leq \bar{G}^{(-1)}((x+\ep)/v(t))-C(\delta)t^{\delta}\)\times\mathbb{P}\(t<\Delta\)\\
\ar\geq\ar\mathbb{E}\[\(1-\bar{F}\(\bar{G}^{(-1)}((x+\ep)/v(t))-C(\delta)t^{\delta}\)\)^{N(t)}\]\times\(1-h(t)\).
\eeqnn
Under condition (\ref{eq:plup}), we continue to derive
\beqnn
\mathbb{P}\(a(v(t))\hat{M}(t)\leq x^{-r_2^{-1}}\)
\ar\geq\ar\mathbb{E}\[\(1-(1+\ep)\bar{G}\(\bar{G}^{(-1)}((x+\ep)/v(t))-C(\delta)t^{\delta}\)\)^{N(t)}\]\times\(1-h(t)\)
\eeqnn
as $t\rightarrow0+$. Similar to (\ref{eq:barF}), we can derive
\beqnn
\lim_{t\rightarrow0+}v(t)\left|{\bar{G}\(\bar{G}^{(-1)}\(x/v(t)\)-{C(\delta)t^{\delta}}\)-x/v(t)}\right|=0
\eeqnn
for any $x>0$, which implies
\beqnn
\bar{G}\(\bar{G}^{(-1)}\(x/v(t)\)-{C(\delta)t^{\delta}}\)\leq(x+\ep)/v(t)
\eeqnn
as $t\rightarrow0+$. Since $\lim_{t\rightarrow0+}h(t)=0$, thus
\beqnn
\mathbb{P}\(a(v(t))\hat{M}(t)\leq x^{-r_2^{-1}}\)
\ar\geq\ar\mathbb{E}\[\(1-(1+\ep) (x+2\ep)/v(t)\)^{N(t)}\](1-\ep)
\eeqnn
as $t\rightarrow0+$.
It follows from Proposition \ref{pro:number} that as $t\rightarrow0+$
\beqnn
\mathbb{P}\(N(t)\leq (1+\ep)v(t)\)\geq 1-O(t^{1-2\al^{*}}).
\eeqnn
Consequently, we continue to get
\beqnn
\mathbb{P}\(a(v(t))\hat{M}(t)\leq x^{-r_2^{-1}}\)\ar\geq\ar\mathbb{E}\[\(1-{(1+\ep) (x+2\ep)}/v(t)\)^{N(t)}1_{\{N(t)\leq (1+\ep)v(t)\}}\](1-\ep)\\
\ar\geq\ar\[\(1-{(1+\ep)(x+2\ep)}/v(t)\)^{(1+\ep)v(t)}\]\times \(1-\ep\)^2
\eeqnn
as $t\rightarrow0+$. Note that
\beqnn
\lim_{t\rightarrow0+}\(1-(1+\ep) (x+2\ep)/v(t)\)^{(1+\ep)v(t)}=e^{-(1+\ep)^2(x+2\ep)}.
\eeqnn
Then the desired result follows by taking $\ep\rightarrow0+$ and applying a change of variable.

(b) follows similarly to the proof of (a).
\end{proof}

As an application of Proposition \ref{pro:ere}, we can show convergence in probability for $\log \hat{M}(t)/{\log v(t)} $ under relaxed condition on asymptotics of $\bar{F}$.
 We introduce the following conditions, which are more general than those in (\ref{eq:plup}) and (\ref{eq:pllow}).
\beqlb\label{eq:4221}
\limsup_{x\rightarrow\infty}\frac{\bar{F}(x)}{R_1(x)}\leq 1
\eeqlb
and
\beqlb\label{eq:4222}
\liminf_{x\rightarrow\infty}\frac{\bar{F}(x)}{R_2(x)}\geq 1,
\eeqlb
where $R_1(x)$ and $R_2(x)$ are regularly varying functions at $\infty$.

\begin{corollary}\label{con_prob}
Suppose that both (\ref{eq:4221}) and (\ref{eq:4222}) hold for regularly varying functions $R_1(x)$ and $R_2(x)$,
each with a common negative index $-r$. Then we have
${\log \hat{M}(t)}/{\log v(t)}\to  {r}^{-1}$ in probability as $t\rightarrow0+$.
\end{corollary}

\begin{proof}
Given any $\ep\in(0,{r}^{-1})$, we choose $\delta$ with $0<\delta<\ep r/(\ep+r^{-1})$, which yields $\ep+{r}^{-1}-(r-\delta)^{-1}>0$. Under condition (\ref{eq:4221}), we have
\beqnn
\bar{F}(x)\leq x^{-r+\delta}
\eeqnn
for $x$ large enough. By Proposition \ref{pro:ere} (\ref{a21}),
\beqnn
\ar\ar\limsup_{t\rightarrow0+}\mathbb{P}\({\log \hat{M}(t)}/{\log v(t)}-{r^{-1}}>\ep\)\\
\ar=\ar{\limsup_{t\rightarrow0+}}\mathbb{P}\(\hat{M}(t)> {v(t)}^{\ep+{r}^{-1}}\)\cr
\ar=\ar1-{\liminf_{t\rightarrow0+}}\mathbb{P}\(v(t)^{-(r-\delta)^{-1}}\hat{M}(t)\leq v(t)^{\ep+{r}^{-1}-(r-\delta)^{-1}}\)\cr
\ar\leq\ar0.
\eeqnn
On the other hand, for the given $\ep\in(0,{r}^{-1})$, we can choose $\delta$ with $0<\delta<\ep r/(r^{-1}-\ep)$, which yields ${r}^{-1}-\ep-(r+\delta)^{-1}<0$. Under condition (\ref{eq:4222}), we have
\beqnn
\bar{F}(x)\geq x^{-r-\delta}
\eeqnn
for $x$ large enough. By Proposition \ref{pro:ere} (\ref{a22}),
\beqnn
\ar\ar\limsup_{t\rightarrow0+}\mathbb{P}\({\log \hat{M}(t)}/{\log v(t)}-{r^{-1}}<-\ep\)\\
\ar=\ar{\limsup_{t\rightarrow0+}}\mathbb{P}\(\hat{M}(t)< {v(t)}^{{r}^{-1}-\ep}\)\cr
\ar=\ar{\limsup_{t\rightarrow0+}}\mathbb{P}\(v(t)^{-(r+\delta)^{-1}}\hat{M}(t)< v(t)^{{r}^{-1}-\ep-(r+\delta)^{-1}}\)\cr
\ar\leq\ar0.
\eeqnn
Therefore,
\beqnn
{\lim_{t\rightarrow0+}}\mathbb{P}\(\big|{\log \hat{M}(t)}/{\log v(t)}-{r^{-1}}\big|>\ep\)=0
\eeqnn
The proof is complete.
\end{proof}

\subsection{Exact speed of CDI in fast regime}\label{fast_regime}
We now present a  condition on $\bar{F}$ so that it is asymptotically exponential function like. Then we show   renormalization results on the  extremal process.  In this subsection we need the condition
\beqlb\label{eq:comuplow}
\lim_{x\rightarrow\infty}\frac{\bar{F}(x)}{ x^{r_3}e^{-r_1x^{r_2}}}= r_4
\eeqlb
for constants $r_1>0$, ${r_2}>0$, $r_3\in\mathbb{R}$ and $r_4>0$.

\begin{assumption}\label{as:v}
Given the constant $\delta$ in \eqref{eq:mod} and $r_2$ in (\ref{eq:comuplow}), the following holds:
\beqlb\label{eq:v}
\(\log v(t)\)^{0\vee (1-r_2^{-1})}t^{\delta}\rightarrow0\,\,\,\,\text{as}\,\,t\rightarrow0+.
\eeqlb
\end{assumption}
\begin{remark}
For $r_2>1$, if $\psi(q)\sim Cq(\log q)^{1+\ep}$ for some $\ep> (1-r_2^{-1})/\delta$, then it follows that $v(t)\sim e^{Ct^{-1/\ep}}$ as $t\rightarrow0+$, which is sufficient for condition (\ref{eq:v}).
\end{remark}

\begin{theorem}\label{th:2}
Under condition (\ref{eq:comuplow}) and Assumption \ref{as:v},
for any $x\in\mathbb{R}$, we have
\beqnn
\lim_{t\rightarrow0+}\mathbb{P}\(r_1\hat{M}(t)^{r_2}-\log v(t)-\log r_4-r_3r_2^{-1}\log\(r_1^{-1}\log v(t)\)\leq x\)=e^{-e^{-x}}.
\eeqnn
\end{theorem}

\begin{proof}
Under condition (\ref{eq:comuplow}), one can show that
\beqnn
\bar{F}^{(-1)}(y)=\(r_1^{-1}[\log (y^{-1})+\log r_4+r_3{r_2^{-1}}\log\(r_1^{-1}\log (y^{-1})\)]\)^{r_2^{-1}}
\eeqnn
is an asymptotic inverse function for $\bar{F}$ as $y\rightarrow0+$.

For any $x>0$, as $t\to 0+$, we have
\beqnn
\bar{F}\(\bar{F}^{(-1)}\(x/v(t)\){\pm}{C(\delta)t^{\delta}}\)\ar\sim\ar r_4\(\bar{F}^{(-1)}\(x/v(t)\){\pm}{C(\delta)t^{\delta}}\)^{r_3}e^{-r_1\(\bar{F}^{(-1)}\(x/v(t)\){\pm}{C(\delta)t^{\delta}}\)^{r_2}}\nonumber\\
\ar\sim\ar r_4\(\bar{F}^{(-1)}\(x/v(t)\)\)^{r_3}e^{-r_1\(\bar{F}^{(-1)}\(x/v(t)\)\)^{r_2}{\mp}r_1r_2\(\bar{F}^{(-1)}\(x/v(t)\)\)^{r_2-1}C(\delta)t^{\delta}}\nonumber\\
\ar\sim\ar x/v(t) \(1{\mp}r_1r_2\(\bar{F}^{(-1)}\(x/v(t)\)\)^{r_2-1}C(\delta)t^{\delta}\)\\
\ar\sim\ar x/v(t) \(1{\mp}r_1r_2\(\log v(t)/r_1\)^{1-r_2^{-1}}C(\delta)t^{\delta}\).
\eeqnn
This fact, combined with Assumption \ref{as:v}, yields
\beqnn
\lim_{t\rightarrow0+}v(t)\left|{\bar{F}\(\bar{F}^{(-1)}\(x/v(t)\){\pm}{C(\delta)t^{\delta}}\)-x/v(t)}\right|=0.
\eeqnn
Thus, condition (\ref{eq:Fa}) holds.

By Theorem \ref{th:uc2}, as $t\rightarrow0+$, we have
\beqlb\label{eq:u2}
\lim_{t\to 0+}\mathbb{P}\Big(\hat{M}(t)\leq \bar{F}^{(-1)}(x/v(t))\Big)= e^{-x}\,\,\,\,\text{for any}\,x>0.
\eeqlb

Given any $\ep>0$, as $t\to 0+$, it is straightforward to verify that
\beqnn
\bar{F}^{(-1)}\Big(\frac{x+\ep}{v(t)}\Big)
\ar\leq \ar\Big(r_1^{-1}\log (v(t)/x)+\frac{\log r_4+r_3{r_2^{-1}}\log\(r_1^{-1}\log (v(t))\)}{r_1}\Big)^{r_2^{-1}}
\leq
\bar{F}^{(-1)}\Big(\frac{x-\ep}{v(t)}\Big).
\eeqnn
It follows from (\ref{eq:u2}) that
\beqnn
e^{-(x+\ep)}\ar=\ar\lim_{t\rightarrow0+}\mathbb{P}\Big(\hat{M}(t)\leq \bar{F}^{(-1)}((x+\ep)/v(t))\Big)\\
\ar\leq\ar\liminf_{t\rightarrow0+}\mathbb{P}\bigg(\hat{M}(t)\leq \Big(r_1^{-1}\log (v(t)/x)+\frac{\log r_4+r_3{r_2^{-1}}\log\(r_1^{-1}\log (v(t))\)}{r_1}\Big)^{r_2^{-1}}\bigg)\\
\ar\leq\ar\limsup_{t\rightarrow0+}\mathbb{P}\bigg(\hat{M}(t)\leq \Big(r_1^{-1}\log (v(t)/x)+\frac{\log r_4+r_3{r_2^{-1}}\log\(r_1^{-1}\log (v(t))\)}{r_1}\Big)^{r_2^{-1}}\bigg)\\
\ar\leq\ar\lim_{t\rightarrow0+}\mathbb{P}\Big(\hat{M}(t)\leq \bar{F}^{(-1)}((x-\ep)/v(t))\Big)=e^{-(x-\ep)}.
\eeqnn
Letting $\ep\to 0+$ we have
\beqnn
e^{-x}\ar=\ar\lim_{t\rightarrow0+}\mathbb{P}\bigg(\hat{M}(t)\leq \Big(r_1^{-1}\log (v(t)/x)+\frac{\log r_4+r_3{r_2^{-1}}\log\(r_1^{-1}\log (v(t))\)}{r_1}\Big)^{r_2^{-1}}\bigg)\\
\ar=\ar \lim_{t\rightarrow0+}\mathbb{P}\Big(r_1\hat{M}(t)^{r_2}-\log v(t)-\log r_4-r_3r_2^{-1}\log\(r_1^{-1}\log v(t)\)\leq -\log x\Big).
\eeqnn
The desired result follows by a change of variable.
\end{proof}

\begin{remark}
To find the asymptotic inverse function for $\bar{F}$ in the proof for Theorem \ref{th:2},  one can use the following heuristic argument.
Assuming that as $x\to\infty$,
\beqnn
y\sim r_4x^{r_3}e^{-r_1x^{r_2}},
\eeqnn
it follows that
\beqlb\label{eq:yfast}
\log y\sim\log r_4+r_3\log x-r_1x^{r_2}
\eeqlb
for $x$ sufficiently large. Consequently, we anticipate
\beqnn
x\sim\(r_1^{-1}\log (y^{-1})\)^{r_2^{-1}}+m(y)
\eeqnn
for $y\to 0+$, where $m(y)=o\big(\(\log (y^{-1})\)^{r_2^{-1}}\big)$ is the adjustment term. Substituting this back into (\ref{eq:yfast}) allows us to solve for $m(y)$ as follows:
\beqnn
\log y\ar\sim\ar \log r_4+r_3\log\Big(\big(r_1^{-1}\log (y^{-1})\big)^{r_2^{-1}}+{m(y)}\Big)-r_1\Big(\big(r_1^{-1}\log (y^{-1})\big)^{r_2^{-1}}+{m(y)}\Big)^{r_2}\\
\ar\sim\ar \log r_4+r_3{r_2^{-1}}\log\(r_1^{-1}\log (y^{-1})\)-r_1\Big(\big(r_1^{-1}\log (y^{-1})\big)^{r_2^{-1}}\Big)^{r_2}-r_1r_2\Big(\big(r_1^{-1}\log (y^{-1})\big)^{r_2^{-1}}\Big)^{r_2-1}m(y),
\eeqnn
which implies that
\beqnn
m(y)\ar\sim\ar\frac{\log r_4+r_3{r_2^{-1}}\log\(r_1^{-1}\log (y^{-1})\)}{r_1r_2\Big(\(r_1^{-1}\log (y^{-1})\)^{1-r_2^{-1}}\Big)}.
\eeqnn
Therefore, we have as $y\to 0+$,
\beqnn
\bar{F}^{(-1)}(y)\ar\sim\ar \(r_1^{-1}\log (y^{-1})\)^{r_2^{-1}}+\frac{\log r_4+r_3{r_2^{-1}}\log\(r_1^{-1}\log (y^{-1})\)}{r_1r_2\(\(r_1^{-1}\log (y^{-1})\)^{1-r_2^{-1}}\)}\\
\ar\sim\ar \(r_1^{-1}\log (y^{-1})+\frac{\log r_4+r_3{r_2^{-1}}\log\(r_1^{-1}\log (y^{-1})\)}{r_1}\)^{r_2^{-1}}.
\eeqnn
  
The above asymptotic inverse can also be derived by Remark \ref{re:1}. Let $z=e^{r_1x^{r_2}}$.
Then
\beqnn
y\ar\sim\ar r_4x^{r_3}e^{-r_1x^{r_2}}
\sim r_4r_1^{-r_3r_2^{-1}}z^{-1}\(\log z\)^{r_3r_2^{-1}}
\eeqnn
as $z\rightarrow\infty$. By Remark \ref{re:1}, we have
\beqnn
z\ar\sim\ar r_4r_1^{-r_3r_2^{-1}}y^{-1}\(\log y^{-1}\)^{r_3r_2^{-1}}.
\eeqnn
Note that $x=\(r_1^{-1}\log z\)^{r_2^{-1}}$. Then
\beqnn
x\ar\sim\ar \(r_1^{-1}\log \(r_4r_1^{-r_3r_2^{-1}}y^{-1}\(\log y^{-1}\)^{r_3r_2^{-1}}\)\)^{r_2^{-1}}.\cr
\eeqnn

\end{remark}

We next present another renormalization equivalent to that in Theorem \ref{th:2}.

\begin{corollary}\label{co:ex}
Under Assumption \ref{as:v}, if
	$\bar{F}(x)\sim {r_4}e^{-r_1x^{r_2}}$ as $x\rightarrow\infty$ for $r_4, r_1, r_2>0$, then for any $x\in\mathbb{R}$,
	\beqnn
	\lim_{t\rightarrow0+}
	\mathbb{P}\big(r_1\hat{M}(t)^{r_2}-\log v(t)-\log {r_4}\leq x\big)=e^{-e^{-x}}.
	\eeqnn
\end{corollary}
\begin{proof}
The asymptotic tail distribution satisfies condition (\ref{eq:comuplow}) by setting $r_3=0$.
Applying Theorem \ref{th:2}, for any $x\in\mathbb{R}$, we derive	
	$$\lim_{t\rightarrow0+}\mathbb{P}\(r_1\hat{M}(t)^{r_2}-\log v(t)-\log {r_4}\leq x\)=e^{-e^{-x}}. $$
The proof is complete.
\end{proof}

For $x>1$ further define functions
\begin{equation}
	\begin{split}
\bar{a}(x)&\equiv\bar{a}(x;r_1,r_2):=r_1{r_2}\({r_1^{-1}}{\log x}\)^{1-{r^{-1}_2}}
    \end{split}
\end{equation}
and
$${\bar{b}}(x)\equiv {\bar{b}}(x;r_1,r_2,r_3,r_4):=r_2 \log x+\log r_4+{r_3}{r_2^{-1}}\log\({r_1^{-1}}{\log x}\).$$

\begin{corollary}\label{le:e}
Under condition (\ref{eq:comuplow}) and Assumption \ref{as:v},
for any $x\in\mathbb{R}$, we have
\beqnn
\lim_{t\rightarrow0+}\mathbb{P}\(\bar{a}(v(t))\hat{M}(t)-{\bar{b}}(v(t))\leq x\)=e^{-e^{-x}}.
\eeqnn
\end{corollary}

\begin{proof}
By Theorem \ref{th:2}, for any $x\in\mathbb{R}$,
we have
\beqnn
\lim_{t\rightarrow0+}\mathbb{P}\(r_1\hat{M}(t)^{r_2}-\log v(t)-\log r_4-r_3r_2^{-1}\log\(r_1^{-1}\log v(t)\)\leq x\)=e^{-e^{-x}},
\eeqnn
which is equivalent to
\beqnn
\lim_{t\rightarrow0+}\mathbb{P}\bigg(\hat{M}(t)\leq {\Big(\frac{x+\log v(t)+\log r_4+r_3r_2^{-1}\log\(r_1^{-1}\log v(t)\)}{r_1}\Big)^{r_2^{-1}}}\bigg)=e^{-e^{-x}}.
\eeqnn
Given any $\ep>0$ and $x\in\mathbb{R}$, as $t\to 0+$, one can see that
\beqnn
\ar\ar\Big(\frac{x-\ep+\log v(t)+\log r_4+r_3r_2^{-1}\log\(r_1^{-1}\log v(t)\)}{r_1}\Big)^{r_2^{-1}}\\
\ar\ar\leq\Big(\frac{\log v(t)}{r_1}\Big)^{r_2^{-1}}\Big(1+\frac{ x+\log r_4+{r_3}{r_2^{-1}}\log\({r_1^{-1}}{\log v(t)}\)}{{r_2}\log v(t)}\Big)\\
\ar\ar=\(\bar{a}(v(t))\)^{-1}\(x+\bar{b}(v(t))\)\\
\ar\ar\leq\Big(\frac{x+\ep+\log v(t)+\log r_4+r_3r_2^{-1}\log\(r_1^{-1}\log v(t)\)}{r_1}\Big)^{r_2^{-1}}.
\eeqnn
Therefore,
\beqnn
e^{-e^{-(x-\ep)}}\ar=\ar\lim_{t\rightarrow0+}\mathbb{P}\bigg(\hat{M}(t)\leq {\Big(\frac{x-\ep+\log v(t)+\log r_4+r_3r_2^{-1}\log\(r_1^{-1}\log v(t)\)}{r_1}\Big)^{r_2^{-1}}}\bigg)\\
\ar\leq\ar\liminf_{t\rightarrow0+}\mathbb{P}\(\hat{M}(t)\leq\(\bar{a}(v(t))\)^{-1}\(x+\bar{b}(v(t))\)\)\\
\ar\leq\ar\limsup_{t\rightarrow0+}\mathbb{P}\(\hat{M}(t)\leq\(\bar{a}(v(t))\)^{-1}\(x+\bar{b}(v(t))\)\)\\
\ar\leq\ar\lim_{t\rightarrow0+}\mathbb{P}\bigg(\hat{M}(t)\leq {\Big(\frac{x+\ep+\log v(t)+\log r_4+r_3r_2^{-1}\log\(r_1^{-1}\log v(t)\)}{r_1}\Big)^{r_2^{-1}}}\bigg)\\
\ar=\ar e^{-e^{-(x+\ep)}}.
\eeqnn
Letting $\ep\rightarrow0+$, the desired result follows.
\end{proof}


\begin{example}
Under Assumption \ref{as:v}, if the initial distribution $F$ follows a standard normal distribution, then
$$\bar{F}(x)\sim {{(2\pi)^{-1/2}}x^{-1}}e^{-{x^2}/{2}}\,\,\,\text{as}\,\,\, x\to\infty.$$
It thus belongs to the fast regime with $r_1=1/2$, $r_2=2$, $r_3=-1$ and
$r_4=1/\sqrt{2\pi}$. Applying Corollary \ref{le:e}, we obtain the re-normalization limit of
$\hat{M}(t)$ with
$$\bar{a}(v(t))=\sqrt{2\log v(t)}\text{\,\, and \,\,} \bar{b}(v(t))=2\log v(t)-\log 2-(\log \pi)/2-(\log\log  v(t))/2.$$
\end{example}

If the equality in condition (\ref{eq:comuplow}) is replaced by inequality, we can similarly derive bounds for the upper and lower limits, respectively, as follows.

\begin{proposition}\label{pro:exre}
Suppose that Assumption \ref{as:v} holds.
\begin{enumerate}[(a)]
\item\label{a51}
If
\beqlb\label{eq:comup}
\limsup_{x\rightarrow\infty}\frac{\bar{F}(x)}{ x^{r_3}e^{-r_1x^{r_2}}}\leq r_4,
\eeqlb
then for any $x\in\mathbb{R}$, we have
\beqnn
\liminf_{t\rightarrow0+}\mathbb{P}\({\bar{a}\(v(t)\)\hat{M}(t)-{\bar{b}}\(v(t)\)}\leq x\){\geq}e^{-e^{-x}}.
\eeqnn
\item\label{a52}
If
\beqlb\label{eq:comlow}
\liminf_{x\rightarrow\infty}\frac{\bar{F}(x)}{ x^{r_3}e^{-r_1x^{r_2}}}\geq r_4,
\eeqlb
then for any $x>0$, we have
\beqnn
\limsup_{t\rightarrow0+}\mathbb{P}\({\bar{a}\(v(t)\)\hat{M}(t)-{\bar{b}}\(v(t)\)}\leq x\){\leq}e^{-e^{-x}}.
\eeqnn
\end{enumerate}
\end{proposition}
The proof of Proposition \ref{pro:exre} follows a similar approach to that of Proposition \ref{pro:ere}. For brevity, we omit the details.

\begin{corollary}\label{cor:1222}
Under Assumption \ref{as:v}, if condition (\ref{eq:comuplow}) holds for ${r_2}=1$,
we have
$$r_1\hat{M}(t)-\log \nu(t)-r_3\log\log \nu(t)-\log \(r_4 r_1^{-r_3}\)\to\xi$$
in distribution as $t\rightarrow0+$,
where $\xi$ follows the Gumbel distribution.
\end{corollary}
\begin{proof}
Notice that $\bar{a}\(v(t)\)=r_1$. Applying Corollary \ref{le:e}, for any $x\in\mathbb{R}$, we have
\beqnn
\ar\ar\lim_{t\rightarrow0+}\mathbb{P}\(r_1\hat{M}(t)-{\bar{b}}\(v(t)\)\leq x\)=\lim_{t\rightarrow0+}\mathbb{P}\( \bar{a}\(v(t)\)\hat{M}(t)-{\bar{b}}\(v(t)\)\leq x\)
=e^{-e^{-x}}.
\eeqnn
\end{proof}

The limiting distributions in Corollary \ref{le:e} and Corollary \ref{cor:1222} depend on the specific form of tail distribution in
(\ref{eq:comuplow}) where the exponential term plays a dominate role in determining the decay rates. It is natural to generalize the power functions to regularly varying functions in the tail distribution and further prove results on convergence in probability. To this end, we introduce the following conditions, 
 which are more general than those in (\ref{eq:comup}) and (\ref{eq:comlow}).

\beqlb\label{eq:tailedexpup}
\limsup_{x\rightarrow\infty}\frac{\bar{F}(x)}{R(x)e^{-r_1x^{r_2}}}\leq 1
\eeqlb
and
\beqlb\label{eq:tailedexplow}
\liminf_{x\rightarrow\infty}\frac{\bar{F}(x)}{R(x)e^{-r_1x^{r_2}}}\geq 1,
\eeqlb
where $r_1>0$, $r_2>0$ are constants and $R(x)$ is a regularly varying function at $\infty$ with index $r_3\in\mathbb{R}$.



	Under relaxed condition on asymptotics of $\bar{F}$ in the fast regime, we can also show convergence in probability results that exhibit a phase transition.

\begin{corollary}\label{cor:12272}
Under Assumption \ref{as:v}, suppose that
both (\ref{eq:tailedexpup}) and (\ref{eq:tailedexplow}) hold for the same constants $r_1,\,r_2$ and possibly different functions $R$.

\begin{enumerate}[(a)]
\item\label{a1} If $r_2>1$, we have $\hat{M}(t)-\({{r_1^{-1}}\log v(t)}\)^{{r_2^{-1}}} \to 0$ in probability as $t\rightarrow0+$.
\item\label{a3} For any $r_2>0$, we have ${\hat{M}(t){\({r_1^{-1}\log v(t)}\)^{-{r_2^{-1}}} }}\to 1$ in probability as $t\rightarrow0+$.
\end{enumerate}
\end{corollary}

\begin{proof}
Suppose that $R=R_1$ with index $\rho_1$ in (\ref{eq:tailedexpup}) and $R=R_2$ with index $\rho_2$ in (\ref{eq:tailedexplow}).
For any $\ep\in(0,1)$ and $\ka>0$, it is easy to verify that by (\ref{eq:tailedexpup}),
$$\bar{F}(x)\leq (1+\ep)R_1(x)e^{-r_1x^{r_2}}\leq x^{\rho_1+\ka}e^{-r_1x^{r_2}}$$
and  by (\ref{eq:tailedexplow}),
$$\bar{F}(x)\geq (1-\ep)R_2(x)e^{-r_1x^{r_2}}\geq x^{\rho_2-\ka}e^{-r_1x^{r_2}}$$
for sufficiently large $x$. Put $\bar{a}(y)\equiv\bar{a}\(y;r_1,r_2\)$, $\bar{b}_1\({y}\)\equiv\bar{b}\({y};r_1,r_2,\rho_1+\ka,1\)$
and $\bar{b}_2\({y}\)\equiv\bar{b}\({y};r_1,r_2,\rho_2-\ka,1\)$.
Let
$$d_1({y})\equiv\[\log r_4+(\rho_1+\ka)r_2^{-1}\log (r_1^{-1}\log {y})\]/[r_1^{r_2^{-1}} r_2(\log {y})^{1-r_2^{-1}}]$$
and
$$d_2({y})\equiv\[\log r_4+(\rho_2-\ka)r_2^{-1}\log (r_1^{-1}\log {y})\]/[r_1^{r_2^{-1}} r_2(\log {y})^{1-r_2^{-1}}].$$

(\ref{a1})
 Observe that $\bar{a}\({v(t)}\)\rightarrow\infty$, ${d}_1({v(t)})\rightarrow0+$ and $\bar{d}_2\({{v(t)}}\)\rightarrow0+$ as $t\rightarrow0+$ provided $r_2>1$.

Under condition (\ref{eq:tailedexpup}), it follows from Proposition \ref{pro:exre} (\ref{a51}) that for any $\ep>0$,
\beqnn
\ar\ar{\limsup_{t\rightarrow0+}}\mathbb{P}\(\hat{M}(t)-\({\log {v(t)}}/{r_1}\)^{{r_2^{-1}}}>\ep\)\cr
\ar\ar={\limsup_{t\rightarrow0+}}\mathbb{P}\(\bar{a}\({v(t)}\)\hat{M}(t)-\bar{b}_1({v(t)})>\bar{a}\({v(t)}\)\(\({\log {v(t)}}/{r_1}\)^{{r_2^{-1}}}+\ep\)-\bar{b}_1({v(t)})\)\cr\\
\ar\ar={\limsup_{t\rightarrow0+}}\mathbb{P}\(\bar{a}\({v(t)}\)\hat{M}(t)-\bar{b}_1({v(t)})>\bar{a}\({v(t)}\)\(- d_1({v(t)})+\ep\)\)\cr
\ar\ar=1-{\liminf_{t\rightarrow0+}}\mathbb{P}\(\bar{a}\({v(t)}\)\hat{M}(t)-\bar{b}_1({v(t)})\leq\bar{a}\({v(t)}\)\(-d_1({v(t)})+\ep\)\)\leq 0,
\eeqnn
Consequently,
\beqlb\label{eq:1191}
\ar\ar\lim_{t\rightarrow0+}\mathbb{P}\(\hat{M}(t)-\({\log {v(t)}}/{r_1}\)^{{r_2^{-1}}}>\ep\)=0.
\eeqlb
Under condition (\ref{eq:tailedexplow}), it follows from Proposition \ref{pro:exre} (\ref{a52}) that for any $\ep>0$,
\beqnn
\ar\ar{\limsup_{t\rightarrow0+}}\mathbb{P}\(\hat{M}(t)-\({\log {v(t)}}/{r_1}\)^{r_2^{-1}}<-\ep\)\cr
\ar\ar\leq{\limsup_{t\rightarrow0+}}\mathbb{P}\(\bar{a}\({v(t)}\)\hat{M}(t)-\bar{b}_2({v(t)})\leq\bar{a}\({v(t)}\)\(\({\log {v(t)}}/{r_1}\)^{r_2^{-1}}-\ep\)- b_2({v(t)})\)\cr
\ar\ar={\limsup_{t\rightarrow0+}}\mathbb{P}\(\bar{a}\({v(t)}\)\hat{M}(t)-\bar{b}_2({v(t)})\leq\bar{a}\({v(t)}\)\(-d_2({v(t)})-\ep\)\)\leq 0.
\eeqnn
Consequently,
\beqlb\label{eq:1192}
\ar\ar{\lim_{t\rightarrow0+}}\mathbb{P}\(\hat{M}(t)-\({\log {v(t)}}/{r_1}\)^{r_2^{-1}}<-\ep\)=0.
\eeqlb
(\ref{eq:1191}) and (\ref{eq:1192}) yield
\beqnn
\ar\ar\lim_{t\rightarrow0+}\mathbb{P}\(\left|\hat{M}(t)-\({\log {v(t)}}/{r_1}\)^{r_2^{-1}}\right|>\ep\)=0.
\eeqnn

(\ref{a3})
Observe that for $r_2>1$, the desired result follows from that of (a).

For $0<r_2\leq1$, we have  as $t\rightarrow0+$
$$\bar{a}\({v(t)}\)\({\log {v(t)}}/{r_1}\)^{{r_2^{-1}}}\rightarrow\infty,$$
$$ d_1\({v(t)}\)=o\(\({\log {v(t)}}/{r_1}\)^{{r_2^{-1}}}\)\,\,\,\text{and}\,\,\, d_2\({v(t)}\)=o\(\({\log {v(t)}}/{r_1}\)^{{r_2^{-1}}}\).$$

Under condition (\ref{eq:tailedexpup}), it follows from Proposition \ref{pro:exre} (\ref{a51}) that for any $\ep>0$,
\beqnn
\ar\ar{\limsup_{t\rightarrow0+}}\mathbb{P}\({\hat{M}(t)}{\({\log {v(t)}}/{r_1}\)^{-{r^{-1}_2}}}-1>\ep\)\cr
\ar\ar=1-{\liminf_{t\rightarrow0+}}\mathbb{P}\({\hat{M}(t)}{\({\log {v(t)}}/{r_1}\)^{-{r^{-1}_2}}}-1\leq\ep\)\cr
\ar\ar=1-{\liminf_{t\rightarrow0+}}\mathbb{P}\(\bar{a}\({v(t)}\)\hat{M}(t)-\bar{b}_1({v(t)})\leq \bar{a}\({v(t)}\)(1+\ep) {\({\log {v(t)}}/{r_1}\)^{{r_2^{-1}}}}-\bar{b}_1({v(t)})\)\cr
\ar\ar=1-{\liminf_{t\rightarrow0+}}\mathbb{P}\(\bar{a}\({v(t)}\)\hat{M}(t)-\bar{b}_1({v(t)})\leq \bar{a}\({v(t)}\)\(\ep {\({\log {v(t)}}/{r_1}\)^{{r_2^{-1}}}}- d_1\({v(t)}\)\)\)\leq0.
\eeqnn
Consequently,
\beqlb\label{eq:1193}
{\lim_{t\rightarrow0+}}\mathbb{P}\({\hat{M}(t)}{\({\log {v(t)}}/{r_1}\)^{-{r^{-1}_2}}}-1>\ep\)=0.
\eeqlb

Under condition (\ref{eq:tailedexplow}), it follows from Proposition \ref{pro:exre} (\ref{a52}) that for any $\ep>0$,
\beqnn
\ar\ar{\limsup_{t\rightarrow0+}}\mathbb{P}\({\hat{M}(t)}{\({\log {v(t)}}/{r_1}\)^{-{r_2^{-1}}}}-1<-\ep\)\cr
\ar\ar\leq{\limsup_{t\rightarrow0+}}\mathbb{P}\(\bar{a}\({v(t)}\)\hat{M}(t)-\bar{b}_2\({v(t)}\)\leq \bar{a}\({v(t)}\)(1-\ep) {\({\log {v(t)}}/{r_1}\)^{{r_2^{-1}}}}-\bar{b}_2\({v(t)}\)\)\cr
\ar\ar={\limsup_{t\rightarrow0+}}\mathbb{P}\(\bar{a}\({v(t)}\)\hat{M}(t)-\bar{b}_2\({v(t)}\)\leq \bar{a}\({v(t)}\)\(-\ep {\({\log {v(t)}}/{r_1}\)^{{r_2^{-1}}}}- d_2\({v(t)}\)\)\)\leq0.
\eeqnn
Consequently,
\beqlb\label{eq:1194}
\ar\ar{\lim_{t\rightarrow0+}}\mathbb{P}\({\hat{M}(t)}{\({\log {v(t)}}/{r_1}\)^{-{r_2^{-1}}}}-1<-\ep\)=0.
\eeqlb
(\ref{eq:1193}) and (\ref{eq:1194}) yield
\beqnn
\lim_{t\rightarrow0+}\mathbb{P}\(\left|{\hat{M}(t)}{\({\log {v(t)}}/{r_1}\)^{-{r_2^{-1}}}}-1\right|>\ep\)=0.
\eeqnn
\end{proof}

\begin{remark}
	 Combing the results in Sections \ref{slow_regime} and  \ref{fast_regime}, we have that the speed of CDI for $\hat{M}(t)$ is of the form $R(v(t))$ where function $R$ is a  regularly varying is the slow regime and a slowly varying function in the fast regime. 
\end{remark}

\subsection{ $\Lambda$-Fleming-Viot process with Brownian spatial motion in $\mathbb{R}^d$}
The modulus of continuity (Proposition \ref{le:weakm}) holds for  the $d$-dimensional $\Lambda$-Fleming-Viot process with underlying Brownian motion.
If the distributions of $X_1(0),X_2(0),\ldots$ are i.i.d. in $\mathbb{R}^d$  with the tail initial distribution $\mathbb{P}\(\|X_1(0)\|\geq C\)$ decaying at different rates, similar convergence results can be established as those derived in Sections {3.1}-{3.4}. We leave the details to the interested reader.

\begin{appendix}
\section{Modulus of continuity and compact support property for $\Lambda$-Fleming-Viot process with Brownian spatial motion}
We first prove a weak version of modulus of continuity result for the $\Lambda$-Fleming-Viot process with Brownian spatial motion in $\mathbb{R}^d$ under Assumption \ref{weaker_con}, and then establish the compact support property.

\begin{proof}[Proof of Proposition \ref{le:weakm}]
 Given $\ep_o>0$,
by (\ref{eq:cdimo}), there exists $s_o>0$ such that for any $0<s<s_o$,
\beqnn
\mathbb{E}\(\sup_{t\in[0,s]}\left|\frac{N(t)}{v(t)}-1\right|^{d_o}\)\leq 1.
\eeqnn
As a result,
\beqnn
\mathbb{P}\(N(s)>v(s)\(v(s)^{\ep_o}+1\)\)\ar\leq\ar \mathbb{P}\(\left|\frac{N(s)}{v(s)}-1\right|>v(s)^{\ep_o}\)\leq v(s)^{-d_o{\ep_o}}.
\eeqnn
For any ${\delta}\in(0,(1-\delta_o)/2)$, it is evident that
\beqnn
\lim_{t\rightarrow0+}\frac{t^{\delta d-d/2-2\delta-1} e^{-t^{-(1-2\delta)}/2}}{e^{-\(t/2\)^{-{\delta_o}}}}=0.
\eeqnn
Consequently, there exists $t_o$ sufficiently small such that for any $0\leq t<t_o$,
\beqlb\label{eq:a3}
{t^{\delta d-d/2-2\delta-1} e^{-t^{-(1-2\delta)}/2}}\leq{e^{-\(t/2\)^{-{\delta_o}}}}.
\eeqlb

Given $n\in[\infty]$, denote by
\beqnn
I_n=\mathbb{P}\(H^{k2^{-n}}\(k2^{-n}-2^{-\ell},k2^{-n}-2^{-\ell-1}\)\geq 2^{-\(\ell+1\)\delta},\,\exists\,\text{some}\,\ell\geq n,\,k\in\[2^{n}\]\).
\eeqnn
For simplicity of notation, we write $v_{\ell}:=v(2^{-\ell})$ and $N_{\ell}:=N\(2^{-n}-2^{-\ell},2^{-n}\)$ for $\ell\geq n$ in the following proof.
Choose $N$ sufficiently large such that $2^{-N}< s_o\wedge t_o$. For any $n\geq N$, we have
\beqnn
I_n\ar\leq\ar 2^{n}\times\sum_{\ell\geq n}\mathbb{P}\(H^{2^{-n}}\(2^{-n}-2^{-\ell},2^{-n}-2^{-\ell-1}\)\geq 2^{-\(\ell+1\)\delta}\)\cr
\ar\leq\ar 2^{n}\times\sum_{\ell\geq n}\mathbb{P}\(H^{2^{-n}}\(2^{-n}-2^{-\ell},2^{-n}-2^{-\ell-1}\)\geq 2^{-\(\ell+1\)\delta},N_{\ell+1}\leq v_{\ell+1}\({ v_{\ell+1}^{\ep_o}}+1\)\)\cr
\ar\ar\qquad+2^{n}\times\sum_{\ell\geq n}\mathbb{P}\(N_{\ell+1}> v_{\ell+1}\({ v_{\ell+1}^{\ep_o}}+1\)\)\cr
\ar\leq\ar 2^{n}\times\sum_{\ell\geq n}2v_{\ell+1}^{ 1+{\ep_o}}\mathbb{P}\(\max_{0\leq s\leq 1}\| B(s)\|\geq 2^{-\(\ell+1\)\(\delta-1/2\)}\)
+2^{n}\times\sum_{\ell\geq n}v_{\ell+1}^{ -d_o\ep_o},
\eeqnn
where $\(B(s)\)_{0\leq s\leq 1}$ is the standard $d$-dimensional Brownian motion with $B(0)=0$. By an estimate on tail  probability of Brownian motion (see \cite[Equation (3.9)]{Da94} or \cite[Lemma 3]{DL98}),  we continue to have
\beqnn
I_n\leq C2^{n}\times\sum_{\ell\geq n}v_{\ell+1}^{ 1+{\ep_o}}2^{ -\(\ell+1\)\(\delta-1/2\)\(d-2\)}e^{-2^{\(\ell+1\)\(1-2\delta\)}/2}+2^{n}\times\sum_{\ell\geq n}v_{\ell+1}^{-d_o\ep_o}.
\eeqnn
Applying Fubini's theorem,
\beqnn
\sum_{n>N}I_n\ar\leq\ar\sum_{n>N}\(C2^{n}\times\sum_{\ell\geq n}v_{\ell+1}^{ 1+{\ep_o}}2^{-\(\ell+1\)\(\delta-1/2\)\(d-2\)}e^{-2^{\(\ell+1\)\(1-2\delta\)}/2}+2^{n}\times\sum_{\ell\geq n}v_{\ell+1}^{-d_o\ep_o}\)\cr
\ar=\ar\sum_{\ell>N}\sum_{N<n\leq\ell}C2^{n}v_{\ell+1}^{ 1+{\ep_o}}2^{-\(\ell+1\)\(\delta-1/2\)\(d-2\)}e^{-2^{\(\ell+1\)\(1-2\delta\)}/2}
+\sum_{\ell>N}\sum_{N<n\leq\ell}2^{n}v_{\ell+1}^{-d_o\ep_o}\cr
\ar\leq\ar\sum_{\ell>N}C2^{\ell+1}v_{\ell+1}^{ 1+{\ep_o}}2^{-\(\ell+1\)\(\delta-1/2\)\(d-2\)}e^{-2^{\(\ell+1\)\(1-2\delta\)}/2}
+\sum_{\ell>N}C2^{\ell+1}v_{\ell+1}^{-d_o\ep_o}\cr
\ar\leq\ar\sum_{\ell>N}Cv_{\ell+1}^{ 1+{\ep_o}}2^{-\(\ell+1\)\(\delta d-d/2-2\delta-1\)}e^{-2^{\(\ell+1\)\(1-2\delta\)}/2}{2^{-(\ell+2)}}
+{\sum_{\ell>N}C2^{\ell+1}v_{\ell+1}^{-d_o\ep_o}}.
\eeqnn
By (\ref{eq:a3}), we have
$${2^{-(\ell+1)\(\delta d-d/2-2\delta-1\)} e^{-2^{(\ell+1)(1-2\delta)}/2}}\leq{e^{-2^{(\ell+2){\delta_o}}}}.$$
Consequently,
\beqnn
\sum_{n>N}I_n\ar\leq\ar\sum_{\ell>N}Cv_{\ell+1}^{ 1+\ep_o}{e^{-2^{(\ell+2){\delta_o}}}}2^{-(\ell+2)}
+\sum_{\ell>N}C2^{\ell+1}v_{\ell+1}^{-d_o\ep_o}.
\eeqnn

Since function $v(t)$ is decreasing and $e^{-t^{-{\delta_o}}}$ is increasing,
(\ref{eq:v1}) implies
\beqlb\label{eq:v3}
\sum_{\ell>N}v_{\ell+1}^{ 1+\ep_o}e^{-2^{\(\ell+2\){\delta_o}}}2^{-\(\ell+2\)}<\infty,
\eeqlb
where we evaluate  $v(t)$  at the right endpoint and evaluate $e^{-t^{-{\delta_o}}}$  at the left endpoint of each subinterval $[2^{-(\ell+2)},2^{-(\ell+1)}]$ for $\ell$ sufficiently large.
Since
$v_{\ell+1}\geq 2\times 2^{\ell+1},$ see e.g. \cite[Proposition 5.1]{LZ3}, then for ${ d_o\ep_o}>1$,
\beqlb\label{eq:v4}
\sum_{\ell>N}2^{\ell+1}v_{\ell+1}^{-d_o\ep_o}\leq 2^{-d_o\ep_o}\sum_{\ell>N} 2^{-\(\ell+1\)\(d_o\ep_o-1\)}<\infty.
\eeqlb
(\ref{eq:v3}) combined with (\ref{eq:v4}) yields
$$\sum_nI_n\leq N+\sum_{n>N}I_n<\infty.$$
By the Borel-Cantelli Lemma, there exists an $\tilde{N}$ large enough so that
for any $n>\tilde{N}$, $\ell\geq n$ and $k\in\[2^{n}\]$,
$$
H^{k2^{-n}}\(k2^{-n}-2^{-\ell},k2^{-n}-2^{-\ell-1}\)< 2^{-\(\ell+1\)\delta}.$$
Consequently, the lookdown representation implies
\beqnn
H^{k2^{-n}}\(\(k-1\)2^{-n},k2^{-n}\)\leq \sum_{\ell\geq n}H^{k2^{-n}}\(k2^{-n}-2^{-\ell},k2^{-n}-2^{-\ell-1}\)
\ar\leq\ar C(\delta) 2^{-n\delta}.
\eeqnn

For any $s,\,t\in [0,1]$ with $0<t-s\leq \Delta:=2^{-\tilde{N}}$, choose $k$ satisfying $2^{-k}\leq t-s<2^{-k+1}$,
and $s_k, t_k\in\{i2^{-k},i\in[2^k]\}$ satisfying $0\leq s_k-s<2^{-k}$, $0\leq t-t_k<2^{-k}$ and $s_k\leq t_k$. Then $t_k-s_k\leq 2^{-k}$.
Let $(s_{\ell})_{\ell\geq k}$ be a decreasing sequence converging to $s$ such that $s_{\ell+1}=s_{\ell}-i_{\ell+1}2^{-(\ell+1)}$ with $i_{\ell+1}\in\{0,1\}$. Let $(t_{\ell})_{\ell\geq k}$ be an increasing sequence converging to $t$ such that $t_{\ell+1}=t_{\ell}+j_{\ell+1}2^{-(\ell+1)}$ for $j_{\ell+1}\in\{0,1\}$.
It follows from the lookdown representation that
\beqnn
H^{s_k}\(s,s_k\)\leq \sum_{\ell\geq k}H^{s_{\ell}}\(s_{\ell+1},s_{\ell}\)\leq \sum_{\ell\geq k}2^{-(\ell+1)\delta}=C(\delta)2^{-k\delta}
\eeqnn
and
\beqnn
H^{t}\(t_k,t\)\leq \sum_{\ell\geq k}H^{t_{\ell+1}}\(t_{\ell},t_{\ell+1}\)\leq \sum_{\ell\geq k}2^{-(\ell+1)\delta}=C(\delta)2^{-k\delta}.
\eeqnn
Therefore,
\beqnn
H^{t}\(s,t\)\leq H^{s_k}\(s,s_k\)+H^{t_k}\(s_k,t_k\)+H^{t}\(t_k,t\)\leq C(\delta)2^{-k\delta}\leq C(\delta)(t-s)^{\delta}.
\eeqnn

Moreover,
we can provide a probability estimate for the random variable $\Delta$ as follows: For any $x>0$,
\beqnn
\mathbb{P}\(\Delta\leq x\)\ar=\ar\mathbb{P}\(2^{-\tilde{N}}\leq x\)=\mathbb{P}\(\tilde{N}\geq \log(1/x)/\log 2\)\leq\sum_{n=\lceil\log(1/x)/\log 2\rceil}^{\infty}I_{n}\nonumber\\
\ar\leq\ar\sum_{\ell=\lceil\log(1/x)/\log 2\rceil}^{\infty}
\[Cv_{\ell+1}^{ 1+\ep_o}{e^{-2^{(\ell+2){\delta_o}}}}2^{-(\ell+2)}+C2^{\ell+1}v_{\ell+1}^{-d_o\ep_o}\]\nonumber\\
\ar\leq\ar
\int_{0}^{2^{-\lceil\log(1/x)/\log 2\rceil}}v(t)^{1+\ep_o}e^{-t^{-\delta_o}}\mathrm{d}t+\frac{2^{-\lceil\log(1/x)/\log 2\rceil\(d_{o}\ep_{o}-1\)}}{1-2^{-(d_o\ep_o-1)}},
\eeqnn
where $\lceil x\rceil$ is the smallest integer greater than or equal to $x$.
\end{proof}

\begin{proof}[Proof of Corollary \ref{co:com}]
The desired result follows from CDI for the $\Lambda$-coalescent and the modulus of continuity from Proposition \ref{le:weakm}. We omit the details.
\end{proof}

\end{appendix}

\noindent{\bf Acknowledgement}
The authors are grateful to Professor Zenghu Li for very helpful discussions.





\end{document}